\documentclass[10pt,leqno]{siamltex}
\usepackage{amsfonts,amsmath,amssymb}
\usepackage{graphics,url,color,epsfig}
\usepackage[hidelinks]{hyperref}
\usepackage{fdsymbol}
\usepackage{stmaryrd}
\usepackage{enumitem}

\usepackage{array}
\usepackage{booktabs}
\usepackage{multirow}

\newtheorem{example}{Example}[section]

\newcommand{\px}[1][x]{\partial_{#1}}
\newcommand{\dx}[1][x]{\,{\rm d}#1}

\newcommand{\mfrac}[1][2]{\frac{1}{2}}

 \allowdisplaybreaks

\title{A stable fast time-stepping method for fractional integral and derivative operators
\thanks{This work was supported by ARC Discovery Project DP150103675}
}

\author{Fanhai Zeng\thanks{School  of  Mathematical   Sciences,
Queensland   University  of  Technology,   Brisbane,   QLD 4001, Australia
(f2.zeng@qut.edu.au).}
\and Ian Turner$^{\dag,}$\thanks{Australian Research Council Centre of Excellence for Mathematical and Statistical Frontiers, Queensland University of
Technology, Brisbane, QLD 4001, Australia (i.turner@qut.edu.au)}
\and Kevin Burrage$^{\dag,}$\thanks{Visiting Professor, Department of  Computer  Science,     University  of
Oxford, OXI 3QD, UK (kevin.burrage@qut.edu.au)}}

\begin{document}

\maketitle

\begin{abstract}
A  unified fast time-stepping method for both fractional integral and
derivative operators is proposed. The fractional operator is decomposed into a local part with
memory length $\Delta T$ and a history part, where the local part is approximated by
the direct convolution method and the history part is approximated by a fast memory-saving method.
The   fast method has $O(n_0+\sum_{\ell}^L{q}_{\alpha}(N_{\ell}))$ active memory and
$O(n_0n_T+ (n_T-n_0)\sum_{\ell}^L{q}_{\alpha}(N_{\ell}))$ operations,
where $L=\log(n_T-n_0)$,
$n_0={\Delta T}/\tau,n_T=T/\tau$, $\tau$ is the stepsize, $T$ is the final time,
and ${q}_{\alpha}{(N_{\ell})}$ is the number of quadrature points used in the
truncated Laguerre--Gauss (LG) quadrature. The error bound of the present fast method is analyzed.
It is shown that the error from the truncated LG quadrature is independent of  the stepsize,
and can be made arbitrarily small
by choosing suitable parameters that are given explicitly.
Numerical examples  are presented  to verify the effectiveness of the current  fast method.
\end{abstract}





\section{Introduction}\label{sec1}
The convolution of the form
\begin{equation}\label{cont-conv}
(k*u)(t)=\int_{0}^tk(t-s)u(s)\dx[s]
\end{equation}
arises in many physical models,
such as  integral equations, integrodifferential equations,  fractional differential equations,
and integer-order differential equations such as  wave propagation
with nonreflecting boundary conditions,
see for example, \cite{BanLopSch17,Diethelm-B10,MetKla00,Pod-B99,WangWJ2017,LubSch02}.

The aim of this paper is to present a stable and fast memory-saving
time-stepping algorithm
for the convolution  \eqref{cont-conv} with  a  kernel $k(t)={t^{\alpha-1}}/{\Gamma(\alpha)}$.
This kind of kernel has found wide applications
in science and engineering   \cite{Diethelm-B10,MetKla00,Pod-B99}.
When $\alpha\geq0$, Eq. \eqref{cont-conv}
gives a fractional
integral of order $\alpha$. If $\alpha<0$, then Eq. \eqref{cont-conv}
can be interpreted
as the Hadamard finite part integral, which is equivalent to the Riemann--Liouville (RL)
fractional derivative of order $-\alpha$, see  {\cite[p. 112]{SamKM-B93}}.

The direct discretization of \eqref{cont-conv} takes the following form
\begin{eqnarray}\label{discrete-conv}
\sum_{k=0}^n\omega_{n,k}u(t_k),{\quad}n=1,2,\ldots,n_T,
\end{eqnarray}
where $\omega_{n,k}$ are the convolution quadrature weights.
The direct computation of \eqref{discrete-conv} requires $O(n_T)$ active memory and
$O(n_T^2)$ operations, which is expensive for long time computations.
The computational difficulty in both memory requirement and computational cost
will increase greatly when
the direct approximation \eqref{discrete-conv} is applied to resolve
high-dimensional time evolution equations and/or a large system of
time-fractional partial differential equations (PDEs) involving
\eqref{cont-conv}, see e.g., \cite{BafHes17b,FordSim2001,LubSch02,YuPK16,ZengZK16}.
However, we believe this computational difficulty for fractional operators
has not been fully addressed in literature.
The short memory principle (see \cite{deng07,Pod-B99})
seems promising to resolve this difficulty, but it has not been widely applied
in  fractional calculus due to its inaccuracy.

Up to now, some progress has been made in  reducing storage requirements
and computational cost for resolving fractional models.
The basic idea is to seek a suitable sum-of-exponentials to
approximate the kernel function $k_{\alpha}(t)=t^{\alpha-1}/\Gamma(\alpha)$, i.e.,
\begin{equation}\label{sum-exp}
k_{\alpha}(t)=t^{\alpha-1}/\Gamma(\alpha)
=\sum_{j=1}^Qw_je^{\lambda_jt} + O(\epsilon),{\quad} t\in[\delta,T]
\end{equation}
where $\delta,T>0$ and $\epsilon>0$ is a given precision.  The key  is to determine
$w_j$ and $\lambda_j$ in \eqref{sum-exp} in order that  the desired accuracy up to
$O(\epsilon)$ can be achieved.

In order to derive \eqref{sum-exp},
Lubich and Sch\"{a}dle \cite{LubSch02} expressed $k_{\alpha}(t)$
in terms of its inverse Laplace transform $k_{\alpha}(t)
= \frac{1}{2\pi i}\int_{\mathcal{C}} \mathcal{L}[k_{\alpha}] e^{t\lambda}\dx[\lambda]$,
where $\mathcal{L}[k_{\alpha}](\lambda)$ denotes the Laplace transform of $k_{\alpha}(t)$ and
$\mathcal{C}$ is a suitable contour.
Then a suitable quadrature was applied to approximate
$\int_{\mathcal{C}}\mathcal{L}[k_{\alpha}]e^{t\lambda}\dx[\lambda]$, which leads to
\eqref{sum-exp}.  This approach can be applied to construct fast methods
for a wide class of nonlocal models.
%
The method in \cite{LubSch02} was then extended to calculate the discrete convolution
\eqref{discrete-conv} in \cite{BanLopSch17,SchLopLub06},
where the coefficients $\omega_{n,j}$ are generated from the generating functions.
A graded mesh version of \cite{LubSch02} was developed in \cite{LopLubSch08} and
an application of \cite{LubSch02} in the simulation of fractional-order viscoelasticity
in complicated arterial geometries was proposed in \cite{YuPK16}.
The storage and computational cost of the fast methods in
\cite{BanLopSch17,LopLubSch08,LubSch02,SchLopLub06} are $O(\log n_T)$ and $O(n_T\log n_T)$,
respectively, which are much less than the direct methods
with   $O(n_T)$ memory and $O(n_T^2)$ operations.
Recently, Baffet and Hesthaven \cite{BaffetH2017,BafHes17b}
have proposed to approximate $\mathcal{L}[k_{\alpha}]$ using
the multipole approximation, which yields \eqref{sum-exp}.


Another approach to derive \eqref{sum-exp} is based on the following integral expression
\begin{equation}\label{kernel}
k_{\alpha}(t)
=\frac{\sin(\alpha\pi)}{\pi}\int_0^{\infty}\lambda^{-\alpha}e^{-t\lambda}\dx[\lambda],
{\quad}\alpha<1.
\end{equation}
Eq. \eqref{kernel} can be derived by inserting $\mathcal{L}[k_{\alpha}]=\lambda^{-\alpha}$  into  Henrici's formula (see  \cite{DAmMurRiz02}), i.e.,
$k_{\alpha}(t) =\frac{1}{2\pi i}\int_{0}^{\infty} \left[\mathcal{L}[k_{\alpha}](\lambda e^{-i\pi})
-\mathcal{L}[k_{\alpha}](\lambda e^{i\pi})\right]e^{-t\lambda}\dx[\lambda]$.
Li \cite{JingLi10} transformed the above integral into
its equivalent form, then   multi-domain Legendre--Gauss quadrature was applied to
approximate  the transformed integral to obtain \eqref{sum-exp}.
Jiang et al  \cite{JiangZZZ17}
combined   Jacobi--Gauss quadrature and   multi-domain Legendre--Gauss quadrature
to discretize \eqref{kernel} for $-1<\alpha<1$, then a global after-processing optimization
technique was applied to further reduce the number of quadrature points,
see also  \cite{YanSunZhang17}.
The exponential sum approximation for $t^{-\beta}(\beta>0)$
has been studied in the literature, which can be used to  design  fast algorithm to approximate
the fractional operators.
The interested readers can refer to \cite{GregoryLucas2010,McLean2016}.
In the references  \cite{BaffetH2017,BafHes17b,BanLopSch17,LopLubSch08,LubSch02,SchLopLub06},
$w_j$ and $\lambda_j$ are complex, however, they are
real in \cite{JingLi10,JiangZZZ17,GregoryLucas2010,McLean2016}.
Apart from the above mentioned fast methods for fractional operators, McLean \cite{McLean12}
proposed to use a degenerate kernel
for evaluating $k_{\alpha}*u(t)$.

In this work, we  derive \eqref{sum-exp} by approximating \eqref{kernel}
using a truncated Laguerre--Gauss (LG) quadrature.
We list  the main contributions of this work as follows.
\begin{itemize}[leftmargin=*]
\item We follow and generalize the framework in  \cite{LubSch02}  to
resolve the short memory principle with a lag-term,
see Section \ref{sec4}.
By choosing suitable  parameters, the current method
can be simplified as that in \cite{BaffetH2017,BafHes17b,JingLi10,JiangZZZ17},
where the time domain is not
divided into exponentially increasing subintervals. This approach simplifies the implementation
of the algorithm.
The present fast method unifies the calculation of the discrete convolutions
to the approximation of both fractional integral and derivative operators
with arbitrary accuracy (see Tables \ref{tb2-5} and  \ref{s5:tb4-2}),
for example, the trapezoidal rule for the
fractional integral operator  (see \cite{Diethelm-B10})
and the L1 method for the fractional derivative operator (see \cite{Diethelm97,SunWu06});
see Figure \ref{eg31fig1}.
  \item
  Given any basis $B$ ($B>1$ is an integer), any  stepsize $\tau$, any  memory length $\Delta T\geq \tau$, and any precisions $\epsilon,\epsilon_0>0$,  the truncation number
  ${q}_{\alpha}(N_{\ell})$ (see \eqref{kappa-N} and \eqref{eq:sec-5-3-5}) of
  the truncated  LG quadrature is determined in order that
  the overall error of the present fast method
 from the truncated LG quadrature is
   $O(\epsilon+\epsilon_0)$, see \eqref{eq:error-3}.
   The truncated LG quadrature and/or a relatively smaller basis $B$
    saves memory and computational cost, see numerical results in
    Tables \ref{tb2}--\ref{tb2-5} and Figure \ref{Fig-2}.
\end{itemize}

We would like to emphasize that  the \emph{truncated} LG quadrature  reduces
the memory and computational cost significantly,
see  Table \ref{tb2}.
The memory and computational cost in \cite{LubSch02,LopLubSch08} can be
halved  due to the symmetry of the trapezoidal rule,
but operations with complex numbers are involved.
In addition, the Gauss--Jacobi and Gauss--Legendre quadrature used
in \cite{JingLi10,JiangZZZ17} may not be truncated.
Furthermore, the discretization error caused from the LG quadrature
is independent of the stepsize and the regularity of the solution
to the considered fractional differential equation (FDE),
and does not appear to be sensitive    to the fractional order
$\alpha\in (-2,1)$ as exhibited in Figure \ref{Fig-2},
which is competitive with the mutlipole approximation \cite{BafHes17b}
in both accuracy and memory requirement, see  Figure \ref{eg32fig1-2}.

In real applications,  analytical solutions to FDEs are unknown
and often non-smooth, and may
have strong singularity at $t=0$, see, for example,
\cite{LiYiChen16,Luchko11,McLMus07,StynesORiGra17}.
In order to resolve the singularity of the solution to the considered FDE,
a graded mesh approach was adopted by some researchers
\cite{LiYiChen16,McLMus07,StynesORiGra17}.
In \cite{StynesORiGra17}, an optimal graded mesh was obtained to achieve
the global convergence  of order $2-\alpha$, $\alpha\in (0,1)$,
which is very effective when
the fractional order is relatively large, but is less effective when
the fractional order tends to zero.
In this paper, we follow Lubich's approach \cite{Lub86}   to deal with
the singularity by introducing  correction terms.

As we mentioned above, a rational approximation was made in \cite{BafHes17b}
to approximate the Laplace transform of the fractional kernel,
while the method in \cite{JingLi10} used the Legendre--Gauss quadrature
to approximate the transformed integral of \eqref{kernel}. These approaches were
initially designed for the fractional integral operator.
However,  we have found that
they can be applied to the RL fractional derivative
operator directly as done in the present work.

This paper is organized as follows.
We follow the approach in \cite{LubSch02} to
present our fast method in Section \ref{sec2}.
The short memory principle  with lag terms is resolved in Section  \ref{sec4},
it unifies the calculation of the discrete convolutions to the approximation of
the fractional integral and derivative operators.  The error analysis of the fast method
is presented in Section \ref{sec:5}, where all the parameters needed in numerical simulations
are explicitly given.
Numerical simulations are presented
to verify the effectiveness of the fast method
in Section \ref{sec:numerical} before the conclusion in the last section.

\section{A stable fast convolution}\label{sec2}
In this section, we follow the approach given in \cite{LubSch02}
to develop our fast convolution.


The goal is to discretize the right-hand side of \eqref{kernel}
using a highly accurate numerical method. It is natural to use   LG quadrature
to approximate \eqref{kernel}, i.e.,
\begin{equation}\label{eq:gauss-2}
k_{\alpha}(t) = \frac{\sin(\alpha\pi)}{\pi}
\int_0^{\infty}\lambda^{-\alpha}e^{-T\lambda}e^{-(t-T)\lambda}\dx[\lambda]
\approx \frac{\sin(\alpha\pi)}{\pi}\sum_{j=0}^N\omega_je^{-(t-T)\lambda_j},
\end{equation}
where $\{\omega_j\}$ and $\{\lambda_j\}$ are
the LG quadrature weights and points that correspond to the weight
function $\lambda^{-\alpha}e^{-T\lambda}$.
The quadrature \eqref{eq:gauss-2}
is exponentially convergent for any $t\geq T$ if $N$ is sufficiently large,
which is discussed in Section \ref{sec:5}.

Denote $t_n=n\tau\,(n=0,1,...,n_T)$ as the grid point, where $\tau$ is the stepsize.
We first restrict ourselves to $\alpha\in[0,1)$.
Using \eqref{eq:gauss-2} and following the idea in \cite{LubSch02},
we present our stable fast  convolution  for approximating
$k_{\alpha}*u(t)$  as  follows:
\begin{itemize}[leftmargin=*]
\item Step 1) Decompose the convolution $k_{\alpha}*u(t)$ as
  \begin{equation}\label{eq:gauss-3-2}\begin{aligned}
k_{\alpha}*u(t)
  =\int_{t-\tau}^tk_{\alpha}(t-s)u(s)\dx[s]+\int_{0}^{t-\tau}k_{\alpha}(t-s)u(s)\dx[s]
  \equiv L^{\alpha}(u,t) + H^{\alpha}(u,t),
  \end{aligned}\end{equation}
  where we call $L^{\alpha}(u,t)$ and $ H^{\alpha}(u,t)$ the local and history
  parts, respectively.  Let $I^{(1)}_{\tau}u(t)$ be the linear interpolation of $u(t)$.
 Then the local  and history parts can be approximated by
  $$L^{\alpha}(u,t_n){\approx} L^{\alpha}(I^{(1)}_{\tau}u,t_n)=L_{\tau}^{(\alpha,n)}u,
  {\quad} H^{\alpha}(u,t_n){\approx}H^{\alpha}(I^{(1)}_{\tau}u,t_n)= H_{\tau}^{(\alpha,n)}u.$$

  \item Step 2)  For every $t=t_n$, let $L$ be the smallest integer satisfying $t_n<2B^L\tau$,
  where $B>1$ is a positive integer.
  For $\ell=1,2,...,L-1$, determine the  integer $q_{\ell}$ such that
  \begin{equation}\label{basis:B}
  s_{\ell}=q_{\ell}B^{\ell}\tau \quad \text{satisfies}\quad
  t_n-s_{\ell}\in[B^{\ell}\tau,(2B^{\ell}-1)\tau].
  \end{equation}
  Set $s_{0}=t_n-\tau$ and $s_L=0$. Then
  $t_{n}-\tau=s_0>s_1>\cdots>s_{L-1}>s_L=0$.
  \item Step 3) Using \eqref{eq:gauss-2}, we approximate  the history part
  $H_{\tau}^{(\alpha,n)}u=H^{\alpha}(I_{\tau}u,t_n)$ by
\begin{equation}\label{eq:gauss-3}\begin{aligned}
H_{\tau}^{(\alpha,n)}u
\approx\frac{\sin(\alpha\pi)}{\pi}\sum_{\ell=1}^L\sum_{j=0}^{N}\omega_j^{(\ell)}
e^{-(t_n-s_{\ell-1}-T_{\ell-1}){\lambda}^{(\ell)}_j}y(s_{\ell-1},s_{\ell},\lambda^{(\ell)}_j)
={}_F{H}^{(\alpha,n)}_{\tau}u
\end{aligned}\end{equation}
with $y(s_{\ell-1},s_{\ell},\lambda)$  given by
\begin{equation}\label{eq:y-1}
y(s_{\ell-1},s_{\ell},\lambda)
=\int_{s_{\ell}}^{s_{\ell-1}}e^{(s-s_{\ell-1})\lambda} I^{(1)}_{\tau} u(s)\dx[s],
\end{equation}
where $\{\omega^{(\ell)}_j\}$ and $\{\lambda^{(\ell)}_j\}$ are
the LG quadrature weights and points that corresponds to
the weight function $\lambda^{-\alpha}e^{-T_{\ell-1}\lambda}$
(see (7.27) in \cite{ShenTW-B11}),  and
$T_{\ell-1}=B^{\ell-1}\tau$   satisfying $t_n-s-T_{\ell-1}\geq 0$
for all $t_n-s\in [B^{\ell-1}\tau,(2B^{\ell}-1)\tau]$, $s\in [s_{\ell},s_{\ell-1}]$.
Here $y(s)=y(s,s_{\ell},\lambda^{(\ell)}_j)$ used in \eqref{eq:gauss-3} that
is defined by  \eqref{eq:y-1} satisfies the following ODE
\begin{equation}\label{ode2}
y'(s)=-\lambda^{(\ell)}_jy(s) + I^{(1)}_{\tau}u(s),{\qquad}y(s_{\ell})=0,
\end{equation}
which
can be exactly solved by the following recursive relation
\begin{equation}\label{ode2-slover}
\begin{aligned}
y(t_{m+1})=&e^{-\lambda^{(\ell)}_j\tau}y(t_m) + e^{-\lambda^{(\ell)}_j\tau}
\int_{t_m}^{t_{m+1}}e^{\lambda^{(\ell)}_j(s-t_m)}I^{(1)}_{\tau} u(s)\dx[s].
\end{aligned}
\end{equation}

\item Step 4)  Calculate the local part
 $L_{\tau}^{(\alpha,n)}u=\int_{t_n-\tau}^{t_n}k_{\alpha}(t_n-s)I^{(1)}_{\tau}u(s)\dx[s]$ with
\begin{equation}\label{interpolation-quadratic}
L_{\tau}^{(\alpha,n)}u = L(I^{(1)}_{\tau}u,t_n)
= \frac{\tau^{\alpha}}{\Gamma(2+\alpha)} (u_{n}-u_{n-1}).
\end{equation}

\end{itemize}
Combining Steps 1)--4), we obtain our fast convolution  for
approximating $k_{\alpha}*u(t)$.
The above fast convolution has the same storage and computational cost as that
in \cite{LubSch02}, the main differences are listed below:
\begin{itemize}
  \item[i)]  The LG quadrature is applied instead of the trapezoidal
  rule to approximate  the history part
  $H^{\alpha}(I_{\tau}u,t)=\int_0^{t-\tau}k_{\alpha}(t-s)I_{\tau}u(s)\dx[s]$, that is,
  the history part $H^{\alpha}(I_{\tau}u,t)$ in  \cite{LubSch02} was approximated by
  \begin{equation}\label{talbot2}\begin{aligned}
  {}_F\hat{H}^{(\alpha,n)}_{\tau}u=
  \mathrm{Im}\left\{\sum_{\ell=1}^L\sum_{j=-N}^{N-1}\hat{\omega}_j^{(\ell)}
  \mathcal{L}[k_{\alpha}](\hat{\lambda}^{(\ell)}_j)
  e^{(t_n-s_{\ell-1}){\hat{\lambda}}^{(\ell)}_j}
  \hat{y}(s_{\ell-1},s_{\ell},\hat{\lambda}^{(\ell)}_j)\right\},
  \end{aligned}\end{equation}
  where $\{\hat{\omega}_j^{(\ell)}\}$ and $\{\hat{\lambda}^{(\ell)}_j\}$ are the weights and quadrature points
  for the Talbot contour $\Gamma_{\ell}$, and
  $\hat{y}(s)=\hat{y}(s,s_{\ell},\hat{\lambda}^{(\ell)}_j)
  =\int_{s_{\ell}}^{s}e^{-(s-s_{\ell-1})\hat{\lambda}^{(\ell)}_j} I^{(1)}_{\tau}u(s)\dx[s]$
   satisfies the following ODE
\begin{equation}\label{ode}
  \hat{y}'(s)=\hat{\lambda}^{(\ell)}_j\hat{y}(s) + I^{(1)}_{\tau}u(s),{\qquad}y(s_{\ell})=0.
\end{equation}
  \item[ii)] We solve a stable ODE \eqref{ode2} instead of a possibly unstable ODE \eqref{ode}
   that may affect the stability and accuracy of \eqref{talbot2}.
   Indeed, for the Talbot contour used in \cite{LubSch02}
(see also the parabolic contour or hyperbolic contour discussed in
\cite{WeidemanTrefethen07}), there exist
$\hat{\lambda}^{(\ell)}_j$'s, whose real parts are positive.
Numerical tests show that \eqref{talbot2} still  works  well since
one may not solve \eqref{ode} for a long time, which reduces the iteration error
from solving \eqref{ode} even though the real part of some
$\hat{\lambda}^{(\ell)}_j$ is positive.
\end{itemize}

The above fast convolution Step 1) -- Step 4) holds for $\alpha<1$,
that is, the RL  fractional derivative
operator of order $-\alpha$ is thus discretized if $\alpha<0$.

\begin{example}\label{eg21}
Let $u(s)=1+s$ in \eqref{eq:gauss-3} and define the relative error
\begin{equation}\label{eq:error1}
e_n = \frac{\big|H^{\alpha}(u,t_n)-{}_F{H}^{(\alpha,n)}_{\tau}u\big|}
{\big|H^{\alpha}(u,t_n)\big|}, {\quad}
\hat{e}_n = \frac{\big|H^{\alpha}(u,t_n)-{}_F\hat{H}^{(\alpha,n)}_{\tau}u\big|}
{\big|H^{\alpha}(u,t_n)\big|}
\end{equation}
where ${}_F{H}^{(\alpha,n)}_{\tau}u$ and ${}_F\hat{H}^{(\alpha,n)}_{\tau}u$ are
defined by  \eqref{eq:gauss-3} and \eqref{talbot2}, respectively.
\end{example}

\begin{figure}[!h]
\begin{center}
\begin{minipage}{0.47\textwidth}\centering
\epsfig{figure=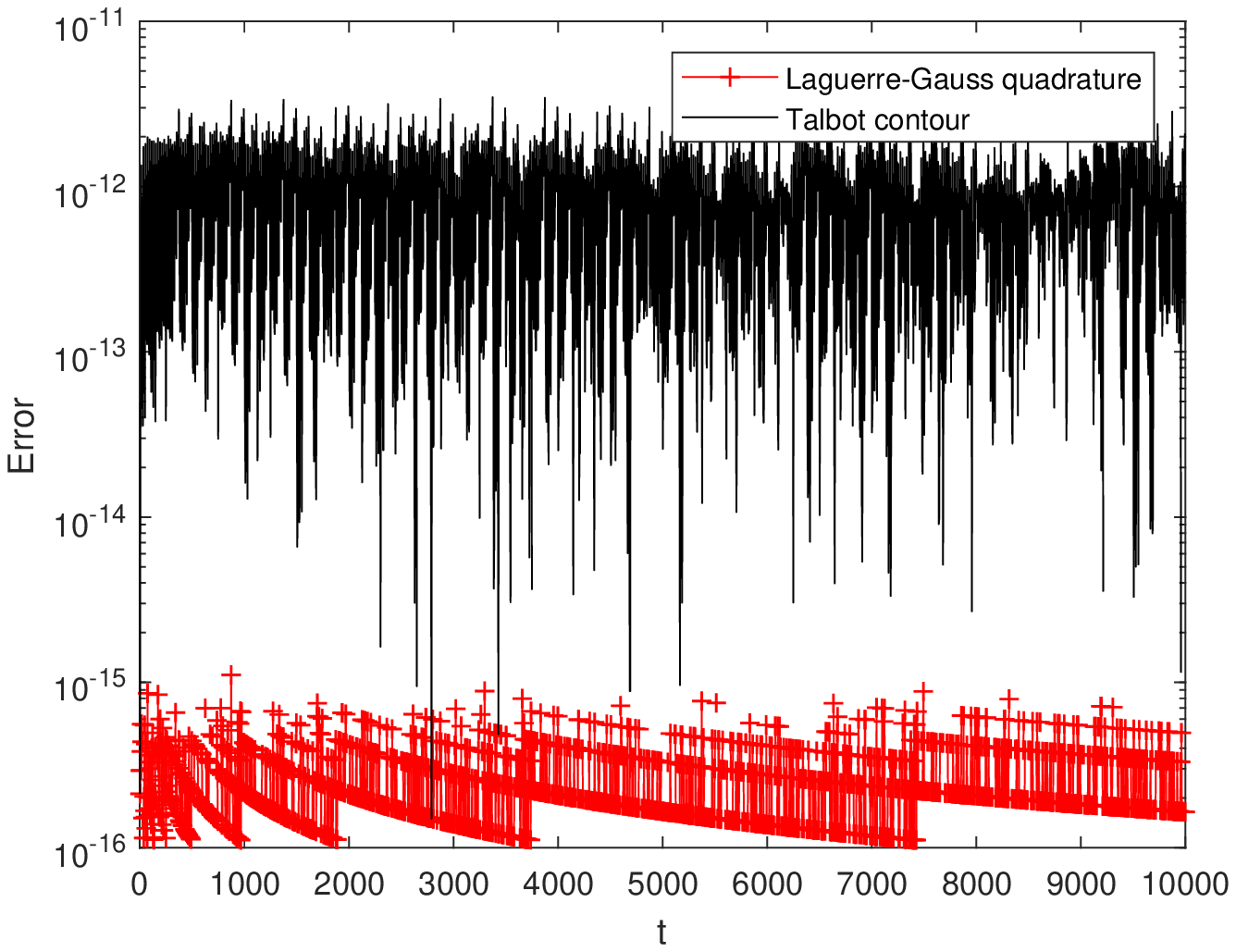,width=5.5cm} \par {(a) $\alpha=-0.5$.}
\end{minipage}
\begin{minipage}{0.47\textwidth}\centering
\epsfig{figure=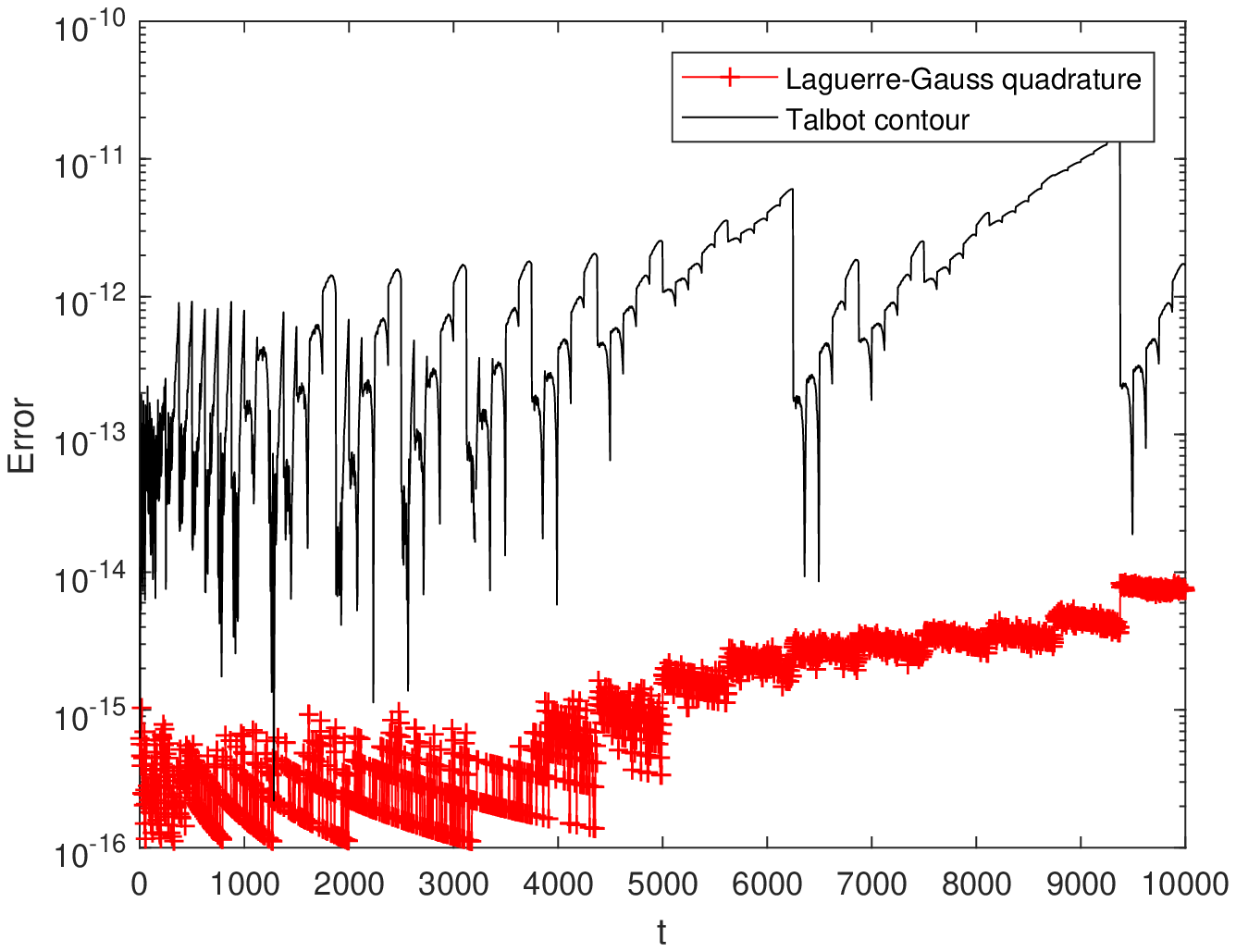,width=5.5cm} \par {(a) $\alpha=0.5$.}
\end{minipage}
\end{center}
\caption{Comparison between the trapezoidal rule based on the
Talbot contour (black curve)
and the LG quadrature (red curve),
$\tau=1$,  $B=5,N=100$. The optimal contour
$z(\theta,t)=t\left(-0.4814 +0.6443(\theta\cot(\theta)+i0.5653\theta)\right)$
obtained in \cite{Weideman06} is applied here, i.e., $\hat{\lambda}_j^{(\ell)}$ in
\eqref{talbot2} is given by $\hat{\lambda}_j^{(\ell)}
= z(\theta_j,N/(2T_{\ell}-\tau))$ with the corresponding weight $\hat{\omega}_j^{(\ell)}
=\px[\theta]z(\theta_j,N/(2T_{\ell}-\tau))$, where $\theta_j=(2j+1)\pi/(2N),j=-N,...,N-1$
and $T_{\ell}=B^{\ell}\tau$, $B=5,N=32$.
\label{eg31fig1}}
\end{figure}

We choose $u(s)=1+s$ in order that the errors $e_n$ and $\hat{e}_n$ mainly come from the
quadrature used in the discretization of the kernel $k_{\alpha}(t)$.
We show the errors  $e_n$ and $\hat{e}_n$ for different  $\alpha$ in Figure \ref{eg31fig1}.
We can see that the LG quadrature shows better accuracy than the trapezoidal
rule based on the Talbot contour in this example.


\section{Short memory principle with lag terms}\label{sec4}
In this section, we  generalize  the fast convolution in the previous section to resolve
the short memory principle (see   \cite{deng07,Pod-B99}).
The error analysis of the method is given in the next section.

Let ${\Delta T}\geq 0$ be a memory length. Divide the convolution
$k_{\alpha}*u(t)$  into the local part
$L^{\alpha}_{{\Delta}T}(u,t)$
and the history part $H^{\alpha}_{{\Delta}T}(u,t)$ as shown below
\begin{eqnarray}
k_{\alpha}*u(t)
&\equiv&L^{\alpha}_{{\Delta}T}(u,t) + H^{\alpha}_{{\Delta}T}(u,t),\label{eq:sec-5-11}\\
L^{\alpha}_{{\Delta}T}(u,t)
&=&\int_{\max\{0,t-{\Delta}T\}}^tk_{\alpha}(t-s)u(s)\dx[s],\label{eq:local}\\
H^{\alpha}_{{\Delta}T}(u,t)
&=&\int_{0}^{\max\{0,t-{\Delta}T\}}k_{\alpha}(t-s)u(s)\dx[s].\label{eq:history}
\end{eqnarray}
If we drop the history part in \eqref{eq:sec-5-11}, then the remaining part
  is the famous short memory principle rule (see \cite{Pod-B99}).
However, the short memory principle
has not been widely applied, since $L^{\alpha}_{{\Delta}T}(u,t)$  is not a good approximation of
$k_{\alpha}*u(t)$. Our goal is to develop a good numerical approximation
of  $H^{\alpha}_{{\Delta}T}(u,t)$, such that  the storage and computational cost
are reduced significantly compared with the direct approximation.
The local  part $L^{\alpha}_{{\Delta}T}(u,t)$
is just the fractional operator defined on $[\max\{0,t-{\Delta}T\},t]$,
which can be discretized directly  by
the known methods, see  \cite{GalGar08,Lub86,LvXu16,SunWu06}.
Next, we introduce how to discretize \eqref{eq:sec-5-11}
efficiently and accurately.

\subsection{Interpolations}\label{sec:3-1}
In this work, the discretization of \eqref{eq:local} and \eqref{eq:history}
is based on the interpolation of $u$.
Specifically,  $L^{\alpha}_{{\Delta}T}(u,t)$ and $H^{\alpha}_{{\Delta}T}(u,t)$
are approximated by $L^{\alpha}_{{\Delta}T}(I^L_{\tau}u,t)$
and $H^{\alpha}_{{\Delta}T}(I^H_{\tau}u,t)$, respectively,
where $I^L_{\tau}$ and $I^H_{\tau}$ are two suitable piecewise interpolation operators,
which will be discussed in the following.


Linear interpolation is simple and has been applied widely in the discretization
of the fractional operators, see \cite{Diethelm-B10,JinRayZhou15,LiZeng2015,SunWu06}.
Several quadratic interpolations have been used to  discretize
the Caputo fractional operators, see, for example, \cite{YanSunZhang17,LiLiangYan17,LvXu16}.
We adopt the quadratic interpolation in \cite{LvXu16}
to illustrate the implementation of the present fast algorithm.
Define the local interpolation operator $\Pi^{j}_{\tau}$ as
\begin{equation}\label{quadratic-0}
\Pi^{j}_{\tau}u(t)= \sum_{k=1}^3u_{j+k-1}F_k^{(j)}(t),{\quad}t\in[t_j,t_{j+1}],j\geq 0,
\end{equation}
where $F^{(j)}_1=\frac{(t-t_{j+1})(t-t_{j+2})}{(t_j-t_{j+1})(t_j-t_{j+2})}$,
$F^{(j)}_2=\frac{(t-t_{j})(t-t_{j+2})}{(t_{j+1}-t_{j})(t_{j+1}-t_{j+2})}$, and
$F^{(j)}_3=\frac{(t-t_{j})(t-t_{j+1})}{(t_{j+2}-t_{j})(t_{j+2}-t_{j+1})}$.
Let $\Pi^{-1}_{\tau}u(t)=\Pi^{0}_{\tau}u(t)$.
Then, for each $n\geq 1$, the quadratic interpolation
$I^{(2,n)}_{\tau}$ is defined  by
\begin{equation}I^{(2,n)}_{\tau}u(t)=\left\{\begin{aligned}\label{quadratic}
&\Pi^{j}_{\tau}u(t), \quad t\in[t_j,t_{j+1}],0 \leq j \leq n-2,\\
&\Pi^{n-2}_{\tau}u(t),{\quad}t\in[t_{n-1},t_{n}].
\end{aligned}\right.\end{equation}
For each $n\geq1$, we define  $I^L_{\tau}$ and $I^H_{\tau}$   as follows
\begin{eqnarray}
I^L_{\tau}u(t)&=&I^{(2,n)}_{\tau}u(t),{\qquad}t\in[\max\{t_{n}-\Delta T,0\},t_n], \label{I-LH-2}\\
I^H_{\tau}u(t)&=&I^{(2,n)}_{\tau}u(t),{\qquad}t\in[0,\max\{t_{n}-\Delta T,0\}]. \label{I-LH-3}
 \end{eqnarray}
Let $j_n = \max\{n-n_0,0\}$.
Then $L^{(\alpha,n)}_{\Delta T,\tau}u=L^{\alpha}_{{\Delta}T}(I^L_{\tau}u,t_n)$
and $H^{(\alpha,n)}_{\Delta T,\tau}u=H^{\alpha}_{{\Delta}T}(I^H_{\tau}u,t_n)$
are given by
\begin{eqnarray}
L^{(\alpha,n)}_{\Delta T,\tau}u
&=&\sum_{j=j_n}^{n-2}\left(b^{(1)}_{n-1-j}u_j
+ b^{(2)}_{n-1-j}u_{j+1}+ b^{(3)}_{n-1-j}u_{j+2}\right)
+ \sum_{j=0}^2d_ju_{n+j-2},\label{discrete-conv-4}\\
H^{(\alpha,n)}_{\Delta T,\tau}u
&=&\sum_{j=0}^{j_n-1}\left(b^{(1)}_{n-1-j}u_j
+ b^{(2)}_{n-1-j}u_{j+1}+ b^{(3)}_{n-1-j}u_{j+2}\right),\label{discrete-conv-5}
\end{eqnarray}
where $b_{n-j-1}^{(k)}=\int_{t_{j}}^{t_{j+1}}k_{\alpha}(t_n-s)F^{(j)}_{k}(s)\dx[s]$
and $d_{j}=\int_{t_{n-1}}^{t_{n}}k_{\alpha}(t_n-s)F^{(n-2)}_{j}(s)\dx[s]$, i.e.,
\begin{eqnarray}
d_0&=&\frac{-\alpha\tau^{\alpha}}{2\Gamma(\alpha+3)},{\quad}
d_1=\frac{\alpha(3+\alpha)\tau^{\alpha}}{\Gamma(\alpha+3)},{\quad}
d_2=\frac{(4+\alpha)\tau^{\alpha}}{2\Gamma(\alpha+3)},\label{discrete-conv-6} \\
b^{(1)}_j&=&\frac{\tau^{\alpha}}{2\Gamma(\alpha)}
\left[\ell_j^{(\alpha+2)} - (2j-1)\ell_j^{(\alpha+1)}
+ j(j-1)\ell_j^{(\alpha)}\right],\label{b1-j}\\
b^{(2)}_j&=&-\frac{\tau^{\alpha}}{\Gamma(\alpha)}
\left[\ell_j^{(\alpha+2)} - 2j\ell_j^{(\alpha+1)}
+ (j+1)(j-1)\ell_j^{(\alpha)}\right],\label{b2-j}\\
b^{(3)}_j&=&\frac{\tau^{\alpha}}{2\Gamma(\alpha)}
\left[\ell_j^{(\alpha+2)} - (2j+1)\ell_j^{(\alpha+1)}
 + j(j+1)\ell_j^{(\alpha)}\right],\label{b3-j}\\
\ell_j^{(\alpha)}&=& \frac{1}{\alpha}\left[(j+1)^{\alpha} - j^{\alpha}\right].
 \end{eqnarray}


Some researchers used other interpolations in the discretization of the fractional
integral and derivative operators, we refer readers to
\cite{LiLiangYan17,Diethelm97,Diethelm-B10,LiZeng2015,FordSim2001,ZayMat2016}.

\subsection{Corrections}
From  \eqref{discrete-conv-4}--\eqref{discrete-conv-5} and for any $\Delta T\geq 0$,
we always have
\begin{equation} \label{discrete-conv-7}
D^{(\alpha,n)}_{\Delta T,\tau}u=L^{(\alpha,n)}_{\Delta T,\tau}u+H^{(\alpha,n)}_{\Delta T,\tau}u
=\sum_{j=0}^{n} w_{n,j}u_{j},
 \end{equation}
where $w_{n,j}$ can be derived from
  \eqref{discrete-conv-6}--\eqref{b3-j}, which do not give  explicitly.
Clearly,  $D^{(\alpha,n)}_{\Delta T,\tau}u$
is just the second-order trapezoidal rule (or the $(2+\alpha)$-order L1 method)
for the fractional integral
of order $\alpha > 0$ (or the RL fractional derivative of order $0<-\alpha < 1$)
if the linear interpolation is used and
$u(t)$ is sufficiently smooth, see \cite{Diethelm97,Diethelm-B10}.
For  quadratic interpolation \eqref{quadratic},
$D^{(\alpha,n)}_{\Delta T,\tau}u$ achieves
the $(3+\alpha)$-order accuracy for  the RL  derivative of order $0<-\alpha<1$
when $u(t)$ is smooth.
However, the approximation $D^{(\alpha,n)}_{\Delta T,\tau}u$ defined by \eqref{discrete-conv-7}
is not a good approximation of $k_{\alpha}*u(t_n)$
when $u(t)$ has strong singularities.

In this work, we follow Lubich's idea (see \cite{Lub86}) to use correction terms to
capture the singularity of the solution $u(t)$ to the considered FODE.
The correction method is based on the
assumption that the solution $u(t)$  has the following form
\begin{equation}\label{solu}
u(t)-u(0)=\sum_{j=1}^{m}c_jt^{\sigma_j} + t^{\sigma_{m+1}}\tilde{u}(t),
{\quad}0<\sigma_j<\sigma_{j+1},
\end{equation}
where $\tilde{u}(t)$ is uniformly bounded for $t\in [0,T]$.
Readers can refer to \cite{Diethelm-B10,Luchko11,JiangLiu-etal12b,Pod-B99}
for more  detailed results of the regularity of FODEs.

Combining \eqref{discrete-conv-7} and \eqref{solu}
gives the following  correction method
\begin{eqnarray}
D^{(\alpha,n,m)}_{\Delta T,\tau}u
&=&D^{(\alpha,n)}_{\Delta T,\tau}u+\tau^{\alpha}\sum_{j=1}^mW_{n,j}(u_j-u_0), \label{direct-conv}
\end{eqnarray}
where
$W_{n,j}$ are the starting weights that are chosen such that
\begin{equation}\label{starting-w}
D^{(\alpha,n)}_{\Delta T,\tau}u+\tau^{\alpha}\sum_{j=1}^mW_{n,j}(u_j-u_0)
=k_{\alpha}*u(t_n)
=\frac{\tau^{\alpha}n^{\sigma_k+\alpha}}{\Gamma(\sigma_k+1+\alpha)},
1\leq k \leq m
\end{equation}
for $u(t)=t^{\sigma_k},0<\sigma_k<\sigma_{k+1}$.  For each $n>0$, one can resolve
$W_{n,j}(1\leq j \leq m)$ from the above linear system and  $W_{n,j}$ are
independent of $\tau$.
The error analysis of the direct method \eqref{direct-conv}
is discussed in Section \ref{sec:5}.
Readers can refer to \cite{DieFFW06,Lub86,ZengZK17}
for more discussions of the correction method.

\subsection{The fast implementation}
In this subsection, we generalize the fast algorithm in Section \ref{sec2} to approximate
$D^{(\alpha,n)}_{\Delta T,\tau}$ defined by \eqref{discrete-conv-7},
which is given as follows.
\begin{itemize}[leftmargin=*]
  \item Step A) Decompose $D^{(\alpha,n)}_{\Delta T,\tau}u$ into two parts
  as  $D^{(\alpha,n)}_{\Delta T,\tau}u
  =L^{(\alpha,n)}_{\Delta T,\tau}u+H^{(\alpha,n)}_{\Delta T,\tau}u$.
  \item Step B) Assume  that ${{\Delta}T}=n_0\tau=t_{n_0}$.
  For  every $t=t_n,n\geq n_0$, let $L$ be the smallest integer satisfying $t_{n-n_0+1}<2B^L\tau$.
  For $\ell=1,2,...,L-1$, determine the  integer $q_{\ell}$ such that
  \begin{equation}\label{L-2}
  s_{\ell}=q_{\ell}B^{\ell}\tau \quad \text{satisfies}\quad
  t_{n-n_0+1}-s_{\ell}\in[B^{\ell}\tau,(2B^{\ell}-1)\tau].
  \end{equation}
  Set $s_{0}=t_{n-n_0+1}-\tau$ and $s_L=0$.
  \item Step C) Let $\hat{t}_n=t_n-{{\Delta}T}+\tau=t_{n-n_0+1}\geq \tau$. Then the history part
  $H^{(\alpha,n)}_{{\Delta}T,\tau}u=H^{\alpha}_{{\Delta}T}(I^H_{\tau}u,t_n)$ is approximated by
  \begin{eqnarray}
H^{(\alpha,n)}_{{\Delta}T,\tau}u
&=&\frac{\sin(\alpha\pi)}{\pi}\sum_{\ell=1}^L\int_0^{\infty}
\lambda^{-\alpha}e^{-(T_{\ell-1}+{{\Delta}T}-\tau)\lambda}e^{-(\hat{t}_n-s_{\ell-1}-T_{\ell-1})\lambda}
y(s_{\ell-1},s_{\ell},\lambda)\dx[\lambda]\nonumber\\
&\approx&\frac{\sin(\alpha\pi)}{\pi}\sum_{\ell=1}^L
\sum_{j=0}^{{q}_{\alpha}{(N_{\ell})}}\omega_j^{(\ell)}
e^{-(\hat{t}_n-s_{\ell-1}-T_{\ell-1}){\lambda}^{(\ell)}_j}
y(s_{\ell-1},s_{\ell},\lambda^{(\ell)}_j)
:={}_FH_{\Delta T,\tau}^{(\alpha,n)}u,\label{eq:sec-5-20}
\end{eqnarray}
where $\{\omega^{(\ell)}_j\}$ and $\{\lambda^{(\ell)}_j\}$ are
the LG quadrature weights and points that correspond  to
the weight function $\lambda^{-\alpha}e^{-(T_{\ell-1}+{{\Delta}T}-\tau)\lambda}$,
${q}_{\alpha}{(N_{\ell})}$ is
defined by \eqref{kappa-N}, and  $y(s_{\ell-1},s_{\ell},\lambda^{(\ell)}_j)$ can be
be obtained exactly  by solving the following linear ODE
\begin{equation}\label{ode-3}
y'(s)=-\lambda^{(\ell)}_jy(s) + I^H_{\tau}u(s),{\quad}y(s_{\ell})=0,
\end{equation}
see also
 \eqref{ode2} and \eqref{ode2-slover}.
  \item Step D) Calculate the local part
  $L^{(\alpha,n)}_{\Delta T,\tau}u=L^{\alpha}_{{\Delta}T}(I^L_{\tau}u,t_n)$.
\end{itemize}

The fast   algorithm for the discretization of
 $k_{\alpha}*u(t)$ is now given by
\begin{equation}\label{fast-conv-correction}
{}_FD^{(\alpha,n,m)}_{\Delta T,\tau}u=L^{(\alpha,n)}_{\Delta T,\tau}u
+{}_FH^{(\alpha,n)}_{\Delta T,\tau}u+\tau^{\alpha}\sum_{j=1}^mW_{n,j}(u_j-u_0),
\end{equation}
where $L^{(\alpha,n)}_{\Delta T,\tau}$  is given by
 \eqref{discrete-conv-4},
${}_FH^{(\alpha,n)}_{\Delta T,\tau}u$ is given by \eqref{eq:sec-5-20},
and the starting weights $W_{n,j}$ are determined by the linear system
\eqref{starting-w}.


Next, we analyze the complexity of the present fast method \eqref{fast-conv-correction}.
For the local part $L^{(\alpha,n)}_{\Delta T,\tau}u$, the memory
requirement is $O(n_0)$ with the computational cost of   $O(n_0(n_T-n_0))$ for all $n_0<n\leq n_T$.
For the history part ${}_FH^{(\alpha,n)}_{\Delta T,\tau}u$,
we have $O(\sum_{\ell}^L{q}_{\alpha}(N_{\ell}))$ active memory and
$O((n_T-n_0)\sum_{\ell}^L{q}_{\alpha}(N_{\ell}))$ operations,
where $L=\log_B{(n_T-n_0)}$.
An additional cost is required to obtain the starting weights $W_{n,j}$
in \eqref{fast-conv-correction}, which can be performed by use of
 fast Fourier transform with  $O(n_T\log(n_T))$ arithmetic operations.
Hence, the overall active memory and computational cost are
$O(n_0+\sum_{\ell}^L{q}_{\alpha}(N_{\ell}))$
and $O(n_0n_T+ (n_T-n_0)\sum_{\ell}^L{q}_{\alpha}(N_{\ell}))$, respectively.


\section{Error analysis} \label{sec:5}
In this section, we analyse the overall discretization error
of the fast method \eqref{fast-conv-correction} in Section   \ref{sec4}. Firstly,
we present the exponential convergence rate of
the LG quadrature used in \eqref{eq:gauss-3} and \eqref{eq:sec-5-20}.
Then we show how to choose  $N_{\ell}$ and ${q}_{\alpha}(N_{\ell})$
used in \eqref{eq:sec-5-20}, such that the desired accuracy
is maintained with the use of the minimum number of the quadrature points.

%

For simplicity, we denote
\begin{equation}\label{eq:sec-5-int}
I^{\alpha}[T,f]=\int_0^{\infty}\lambda^{\alpha}e^{-T\lambda}f(\lambda)\dx[\lambda]{\quad}
\text{and}{\quad}
I^{\alpha}[f]=\int_0^{\infty}\lambda^{\alpha}e^{-\lambda}f(\lambda)\dx[\lambda].
\end{equation}
The LG quadrature for $I^{\alpha}[T,f]$ and $I^{\alpha}[f]$
are given by (see \cite{ShenTW-B11})
\begin{equation}\label{eq:sec-5-GL}
Q_N^{\alpha}[T,f]=T^{-\alpha-1}\sum_{j=0}^N{w}^{(\alpha)}_jf({\lambda}_j/T),{\quad}
Q_N^{\alpha}[f]=\sum_{j=0}^N{w}^{(\alpha)}_jf({\lambda}_j),
\end{equation}
where $\lambda_j$ are the roots of the Laguerre polynomial
$L_{N+1}^{(\alpha)}({\lambda})$, and  ${w}^{(\alpha)}_j$ are the corresponding weights given by
\begin{eqnarray}
{w}^{(\alpha)}_j=\frac{\Gamma(N+\alpha+1){\lambda}_j}{(N+\alpha+1)(N+1)!}
\big(L_N^{(\alpha)}({\lambda}_j)\big)^{-2}.\label{Laguerre-weghts}
\end{eqnarray}

We show   the convergence of  the quadrature
$Q_N^{\alpha}[T,e^{-t\lambda}]$ and the property of the
quadrature weight ${w}^{(\alpha)}_j$ (see \eqref{Laguerre-weghts})
in the following two theorems. The proofs are given in Appendix \ref{appendix-A}.
%
\begin{theorem}\label{thm-1}
Let $t\geq 0,T>0$, $\alpha>-1$, and $N$ be sufficiently large. Then
\begin{equation}\label{eq:sec-5-2-2}
\left|I^{\alpha}[T,e^{-t\lambda}]- Q_N^{\alpha}[T,e^{-t\lambda}]\right|
\leq C_{\alpha,N}T^{-\alpha-1}\left(\frac{t/T}{1+t/T}\right)^{2N},
\end{equation}
where $C_{\alpha,N}$ is bounded $-1<\alpha\leq 0$
and $C_{\alpha,N}\leq C_{\alpha} N^{\alpha}$ for $\alpha>0$.
\end{theorem}
%
\begin{theorem}\label{thm-2}
Let  $\alpha>-1$ and $w^{(\alpha)}_j$ be defined by \eqref{Laguerre-weghts}.
If $N$ and $j$ are sufficiently large,  then there exists a positive constant
$C$ independent of $N$ such that
\begin{equation}\label{L-wj}
{w}^{(\alpha)}_j\leq C(N+1)^{\alpha} e^{-0.25\pi^2(j+1)^2/(N+1)}.
\end{equation}
\end{theorem}

Given a precision $\epsilon_0>0$, the truncated LG  quadrature  is given by
\begin{equation}\label{eq:sec-5-3}
Q_{N,\epsilon_0}^{\alpha}[T,f]=T^{-\alpha-1}
\sum_{j=0}^{q_{-\alpha}(N)}{w}^{(\alpha)}_jf({\lambda}_j/T),
\end{equation}
where $q_{-\alpha}(N)$   is a positive number given by
\begin{equation}\label{kappa-N}
q_{-\alpha}(N)
=\min\Big\{N,\left\lceil{2}{\pi^{-1}}\sqrt{(N+1)\log((N+1)^{\alpha}\epsilon^{-1}_0)}\right\rceil-1\Big\}.
\end{equation}
From \eqref{L-wj}, we can find the smallest integer $j$ satisfying
$(N+1)^{\alpha} e^{-0.25\pi^2(j+1)^2/(N+1)}\leq \epsilon_0$,
which yields \eqref{kappa-N}.
We  choose $\epsilon_0=10^{-16}$ in this paper.

\begin{figure}[!h]
\begin{center}
\begin{minipage}{0.47\textwidth}\centering
\epsfig{figure=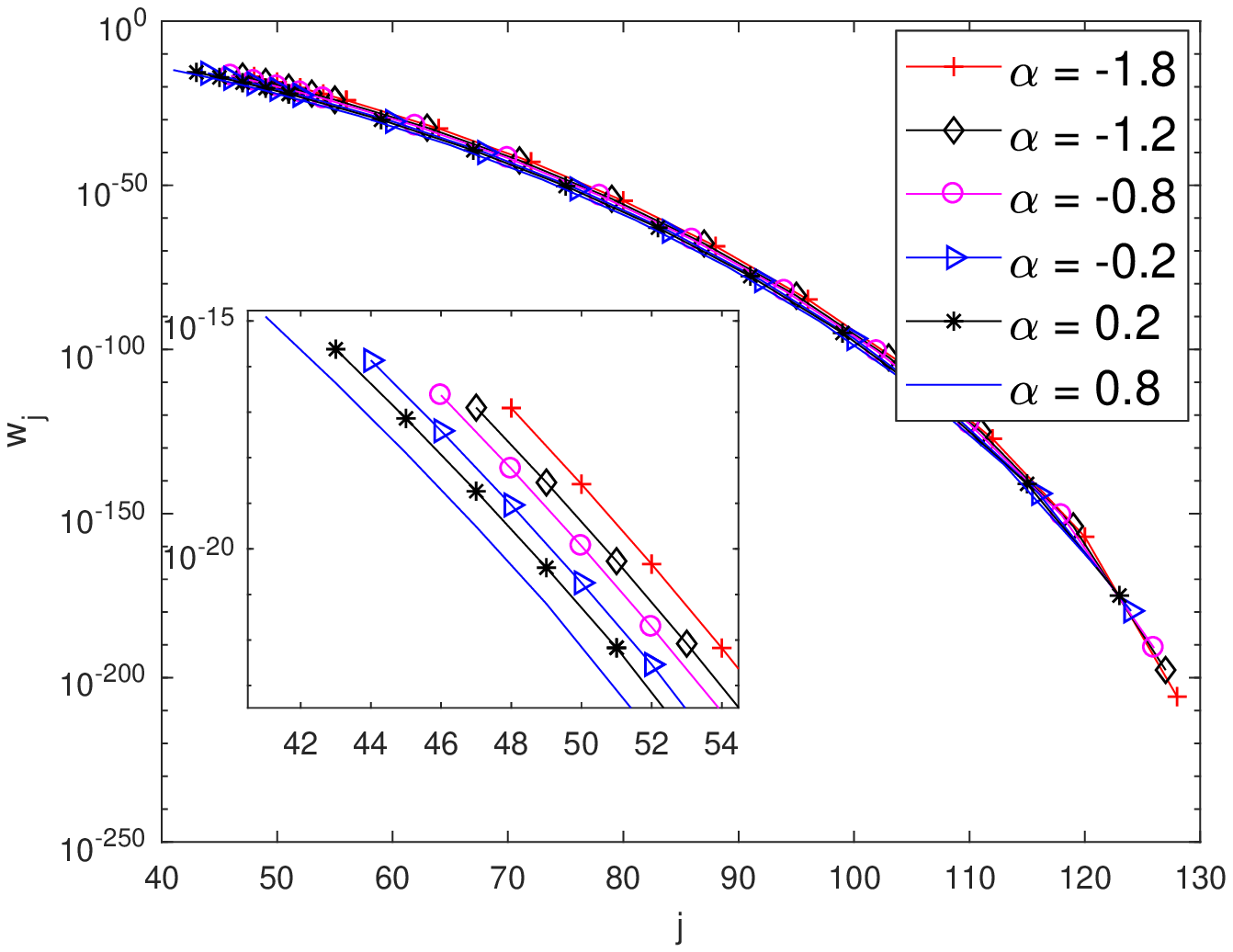,width=5.5cm}
{\par (a) $N=128$.}
\end{minipage}
\begin{minipage}{0.47\textwidth}\centering
\epsfig{figure=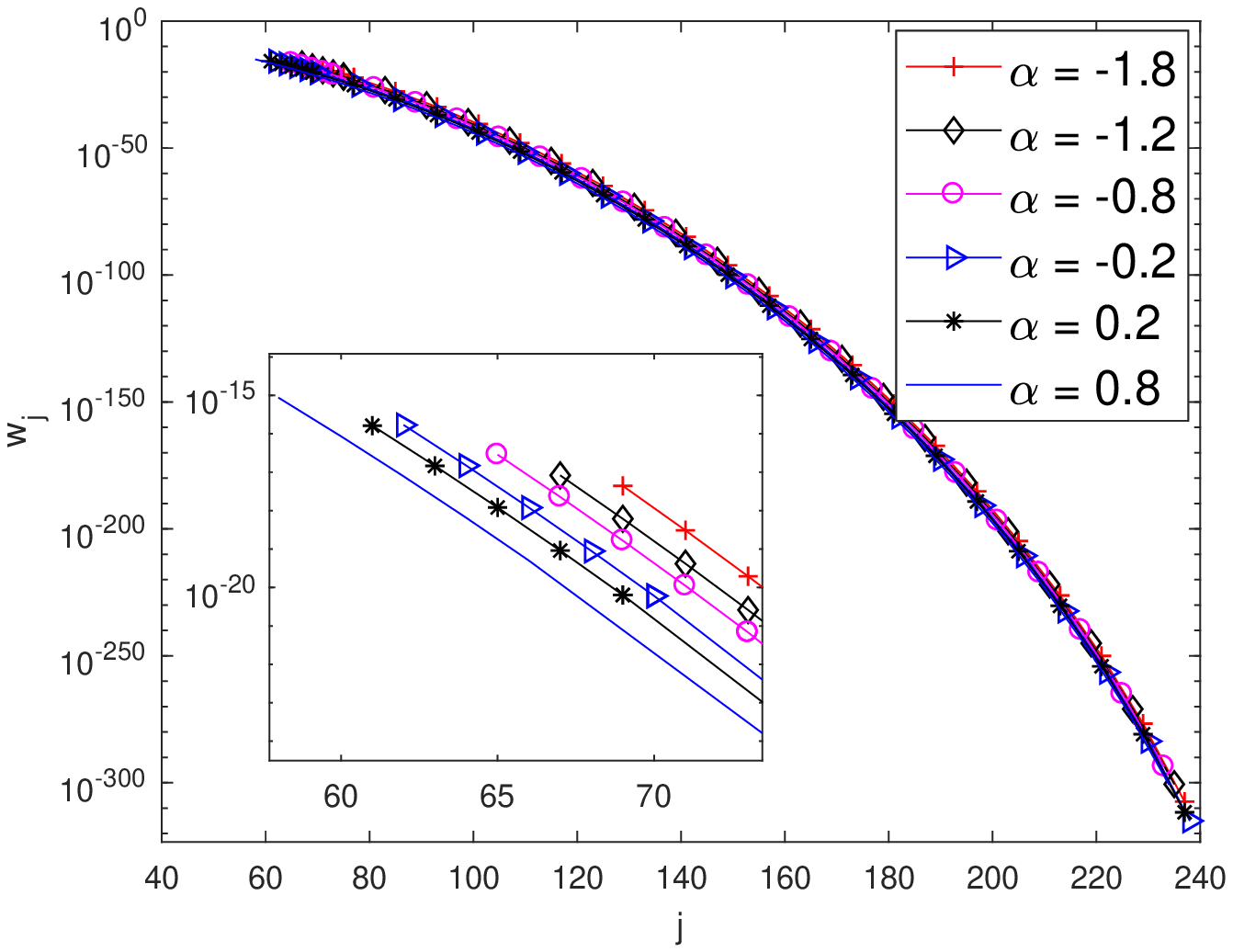,width=5.5cm}
{\par (b) $N=256$.}
\end{minipage}
\end{center}
\caption{The  exponential decay of the   quadrature weights ${w}^{(-\alpha)}_j$
defined by \eqref{Laguerre-weghts}.  \label{Fig-1}}
\end{figure}

Figure \ref{Fig-1} (a) shows the exponential decay of
${w}^{(-\alpha)}_j$ when $j\geq{q}_{\alpha}(128)$
for different fractional orders $\alpha=-1.8,-1.2,-0.8,-0.2,0.2,0.8$.
Figure \ref{Fig-1} (b) displays similar behaviors as shown in Figure \ref{Fig-1} (a).
We can see that Eq. \eqref{kappa-N} works well, and
${q}_{\alpha}(N)$ is not very sensitive to the fractional order $\alpha\in (-2,1)$.
For example, ${q}_{\alpha}(128)=(48,47,46,44,43,41)$ (or ${q}_{\alpha}(256)=(69,67,65,62,61,58)$)
for $\alpha=(-1.8,-1.2,-0.8,-0.2,0.2,0.8)$,
and ${q}_{\alpha}(N) \ll N$ when $N$ is sufficiently large.
For $\alpha \leq -2$, similar results are obtained, which is not shown here.

For notational simplicity, we denote
$$\widehat{T}_{\ell}=T_{\ell-1}+{\Delta} T-\tau,{\quad}
\widehat{H}^n_{\ell}(s)=\hat{t}_n-T_{\ell-1}-s.$$
Next, we investigate how to estimate $N_{\ell}$  in  \eqref{eq:sec-5-20},
such that the   LG quadrature
$Q_{N_{\ell}}^{-\alpha}[\widehat{T}_{\ell},e^{-\widehat{H}^n_{\ell}(s)\lambda}]$
to  the integral
$I^{-\alpha}[\widehat{T}_{\ell},e^{-\widehat{H}^n_{\ell}(s)\lambda}]$
preserves the accuracy up to $O(\epsilon)$ for all $s\in [s_{\ell},s_{\ell-1}]$.

From Theorem \ref{thm-1},
we know that the error of $Q_{N_{\ell}}^{-\alpha}[\widehat{T}_{\ell},e^{-\widehat{H}^n_{\ell}(s)\lambda}]$
 mainly depends on
the following term
$$\left({\widehat{H}^n_{\ell}(s)}/{\widehat{T}_{\ell}}\right)^{2N}
=\left(\frac{\hat{t}_n-s-T_{\ell-1}}{T_{\ell-1}+{\Delta} T-\tau}\right)^{2N},
\quad s\in [s_{\ell},s_{\ell-1}].$$
Using \eqref{L-2} and $T_{\ell-1}=B^{\ell-1}\tau$ gives
\begin{equation}\label{eq:sec-5-3-3-0}
0\leq \frac{\hat{t}_n-s-T_{\ell-1}}{T_{\ell-1}+{\Delta} T-\tau}\leq
\frac{2B-1-B^{1-\ell}}{1+B^{1-\ell}({\Delta}T/\tau-1)}
=\mathcal{T}_{\ell} \leq 2B-1,{\quad} \forall\ell\geq 1.
\end{equation}
Using the above inequality and Eq. \eqref{eq:sec-5-2-2}  yields
\begin{equation}\label{eq:sec-5-3-6}\begin{aligned}
&\big|I^{-\alpha}[\widehat{T}_{\ell},e^{-\widehat{H}^n_{\ell}(s)\lambda}]
-Q_{N_{\ell}}^{-\alpha}[\widehat{T}_{\ell},e^{-\widehat{H}^n_{\ell}(s)\lambda}]\big|\\
\leq& C_{\alpha,N_{\ell}}\widehat{T}_{\ell}^{\alpha-1}
\left(\frac{{\widehat{H}^n_{\ell}(s)}/{\widehat{T}_{\ell}}}
{1+{\widehat{H}^n_{\ell}(s)}/{\widehat{T}_{\ell}}}\right)^{2N_{\ell}}
\leq   C_{\alpha,N_{\ell}}\widehat{T}_{\ell}^{\alpha-1}
\left(\frac{\mathcal{T}_{\ell}}{1+\mathcal{T}_{\ell}}\right)^{2N_{\ell}}.
\end{aligned}\end{equation}
Since the relative error of \eqref{eq:sec-5-3-6} is independent of $\widehat{T}_{\ell}$,
we can let $(\mathcal{T}_{\ell}/(\mathcal{T}_{\ell}+1))^{2N_{\ell}}\leq\epsilon$,
which yields the minimum   $N_{\ell}$ given by
\begin{equation}\label{eq:sec-5-3-5}
N_{\ell} = \left\lceil\frac{\log\epsilon}
{2\log(\frac{\mathcal{T}_{\ell}}{\mathcal{T}_{\ell}+1})}\right\rceil,
\quad \mathcal{T}_{\ell}=\frac{2B-1-B^{1-\ell}}{1+B^{1-\ell}({\Delta}T/\tau-1)}.
\end{equation}

From \eqref{eq:sec-5-3} and \eqref{eq:sec-5-3-6}, we derive that the pointwise error of the
truncated quadrature
$Q_{N_{\ell},\epsilon_0}^{-\alpha}[\widehat{T}_{\ell},e^{-\widehat{H}^n_{\ell}(s)\lambda}]$
for all $s\in[s_{\ell},s_{\ell-1}]$ is given by
\begin{equation}\label{eq:sec-5-3-7}\begin{aligned}
E^{-\alpha,\ell}_{\epsilon,\epsilon_0}[e^{-\widehat{H}^n_{\ell}(s)\lambda}]
=I^{-\alpha}[\widehat{T}_{\ell},e^{-\widehat{H}^n_{\ell}(s)\lambda}]
-Q_{N_{\ell},\epsilon_0}^{-\alpha}[\widehat{T}_{\ell},e^{-\widehat{H}^n_{\ell}(s)\lambda}]
=O(\epsilon) + O(\epsilon_0),
\end{aligned}\end{equation}
where $Q_{N_{\ell},\epsilon_0}^{-\alpha}$ is defined by \eqref{eq:sec-5-3} and
$N_{\ell}$ is given by \eqref{eq:sec-5-3-5}.

Next, we present the  error bound of the fast method  in Section \ref{sec4}.
Denote by
\begin{equation}\label{eq:error-1}\begin{aligned}
R^{(n)}=\int_{0}^{t_n}k_{\alpha}(t_n-s)u(s)\dx[s]-D^{(\alpha,n,m)}_{\Delta T,\tau}u.
\end{aligned}\end{equation}
Note that the above discretization error $R^{(n)}$
depends on the smoothness of $u(t)$ and the  discretization method
$D^{(\alpha,n,m)}_{\Delta T,\tau}u$. If $u(t)$ is sufficiently smooth, no correction
terms are needed to achieve $(2+\alpha)$-order (or $(3+\alpha)$-order)
accuracy if  linear (or quadratic) interpolation is applied for $\alpha<0$.
If $u(t)$ satisfies \eqref{solu}, then the global
$(2+\alpha)$-order (or $(3+\alpha)$-order) accuracy can be achieved
for $\sigma_{m+1}\geq 2$ (or $\sigma_{m+1}\geq 3$).
In numerical simulations, the condition $\sigma_{m+1}\geq 2$ (or $\sigma_{m+1}\geq 3$)
does not need to be satisfied,  $(2+\alpha)$-order (or $(3+\alpha)$-order)
accurate numerical solutions are observed far from $t=0$.
In the numerical simulations, only a small number of correction terms are
sufficient to achieve very accurate numerical solutions;
see numerical simulations in the following section
and see also related results in \cite{ZengZK17}.


From \eqref{eq:sec-5-20} and \eqref{eq:sec-5-3-7}, we have
\begin{equation}\label{eq:error-2}\begin{aligned}
\big|H^{\alpha}_{{\Delta}T}&(I^H_{\tau}u,t_n)-{}_FH_{\Delta T,\tau}^{(\alpha,n)}u\big|
=\Big|\frac{\sin(\alpha\pi)}{\pi}\sum_{\ell=1}^L
E^{-\alpha,\ell}_{\epsilon,\epsilon_0}
[e^{-(\hat{t}_n-T_{\ell-1}-s_{\ell-1})\lambda}y(s_{\ell-1},s_{\ell},\lambda)]\Big|\\
=& \Big|\frac{\sin(\alpha\pi)}{\pi}\sum_{\ell=1}^L
\int_{s_{\ell}}^{s_{\ell-1}}
E^{-\alpha,\ell}_{\epsilon,\epsilon_0}[e^{-\widehat{H}^n_{\ell}(s)\lambda}]
I_{\tau}^{H}u(s)\dx[s]\Big|\\
\leq& C(\epsilon+\epsilon_0)\int_{0}^{\hat{t}_n}\big|I_{\tau}^{H}u(s)\big|\dx[s]
\leq C\max\{0,t_{n+1}-\Delta T\} \|u\|_{\infty}(\epsilon+\epsilon_0).
\end{aligned}\end{equation}
Combining \eqref{eq:error-1} and \eqref{eq:error-2} yields
\begin{equation}\label{eq:error-3}\begin{aligned}
\Big|k_{\alpha}*u(t_n)-{}_FD_{\Delta T,\tau}^{(\alpha,n,m)}u\Big|
=&\Big|H^{\alpha}_{{\Delta}T}(I^H_{\tau}u,t_n)-{}_FH_{\Delta T,\tau}^{(\alpha,n)}u
+ R^{(n)}\Big|\\
\leq& C\max\{0,t_{n+1}-\Delta T\}\|u\|_{\infty} (\epsilon+\epsilon_0) + \big|R^{(n)}\big|,
\end{aligned}\end{equation}
where the error in \eqref{eq:error-3} originates from two parts: the LG quadrature
for discretizing $I^{-\alpha}[\widehat{T}_{\ell},e^{-\widehat{H}^n_{\ell}(s)\lambda}]$
(see \eqref{eq:sec-5-3-7}) and the discretization error defined by \eqref{eq:error-1}.

Next, we numerically study the error caused by the LG quadrature.
Let $m=0$ and $u(t)=1+t$. Then  $R^{(n)}$ in \eqref{eq:error-3} is zero.
Denote the relative pointwise error
$$e_n=\left(k_{\alpha}*u(t_n)\right)^{-1}\Big|k_{\alpha}*u(t_n)
-{}_FD_{\Delta T,\tau}^{(\alpha,n,0)}u\Big|,{\quad}1\leq n \leq n_T=T/\tau,$$
where $k_{\alpha}*u(t_n)=t_n^{\alpha}/\Gamma(1+\alpha)+t_n^{\alpha+1}/\Gamma(2+\alpha)$
for $u(t)=1+t$.

Given a precision $\epsilon= 10^{-10}$,
the maximum relative error
$\|e\|_{\infty}=\max_{0\leq n\leq n_T}\big|e_n\big|$,
the total number of the quadrature points $\sum {N_{\ell}}$,
and the total number of the truncated quadrature points  $\sum{q}_{\alpha}(N_{\ell})$
are shown in Table \ref{tb2} for different basis $B$,
$\alpha=-0.5,\tau=0.1,{\Delta}T = 1$, and $T=10^4$.
We can see that the truncated LG quadrature saves memory. A relatively smaller basis $B$ needs less
quadrature points and thus saves memory, which can be explained from
\eqref{eq:sec-5-3-3-0} and \eqref{eq:sec-5-3-6}. Eq. \eqref{eq:sec-5-3-6}
implies a faster convergence as $B$ decreases, due to
$T_{\ell}\approx 2B-1$, $\ell$ is  sufficiently large. Hence, a relatively smaller $B$
means that  smaller $N_{\ell}$ are needed to achieve high accuracy
that leads to the use of less LG quadrature points.

We change the precision $\epsilon$ and show
the corresponding relative maximum errors $\|e\|_{\infty}$ in Table \ref{tb2-5}.
We can see that    $\|e\|_{\infty}$  increases as $\epsilon$
increases and the total number of the quadrature points are reduced.
It also shows that much better results are obtained than the theoretical prediction,
see also the related results in Table \ref{s5:tb4-2}.

\begin{table}[!h]
\caption{The maximum relative error
$\|e\|_{\infty}$ under the precision $\epsilon= 10^{-10}$,
$\alpha=-0.5,\tau=0.1,{\Delta}T = 1$, and $T=10^4$.}\label{tb2}
\centering\footnotesize
\begin{tabular}{|c|c|c|c|c|c|c|c|c|c|}
\hline
$B$&$\sum N_{\ell}$ &$\sum {q}_{\alpha}(N_{\ell})$ & $\|e\|_{\infty}$
&$B$&$\sum N_{\ell}$ &$\sum {q}_{\alpha}(N_{\ell})$ & $\|e\|_{\infty}$  \\\hline
2 &583    &389 &6.8651e-13 & 30 &3333    &545 &7.2134e-13\\
3 &605    &323 &6.6718e-13 & 40 &3583    &514 &7.3190e-13\\
4 &687    &318 &7.4246e-13 & 50 &4521    &576 &7.3391e-13\\
5 &783    &320 &7.4754e-13 & 60 &5463    &637 &7.4632e-13\\
8 &1151   &370 &7.4377e-13 & 70 &6406    &689 &7.3451e-13\\
10&1246   &359 &7.1388e-13 & 80 &7345    &737 &7.4678e-13\\
15&1598   &376 &6.9954e-13 & 90 &8288    &785 &7.3892e-13\\
20&2173   &438 &7.0761e-13 & 100&9231    &829 &7.4941e-13\\ \hline
\end{tabular}
\end{table}


\begin{figure}[!h]
\begin{center}
\begin{minipage}{0.47\textwidth}\centering
\epsfig{figure=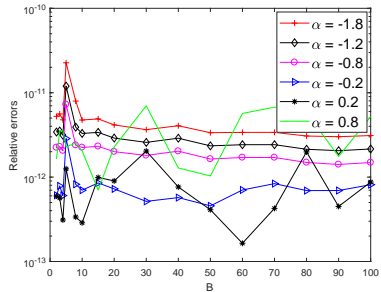,width=5.5cm}
{\par (a) }
\end{minipage}
\begin{minipage}{0.47\textwidth}\centering
\epsfig{figure=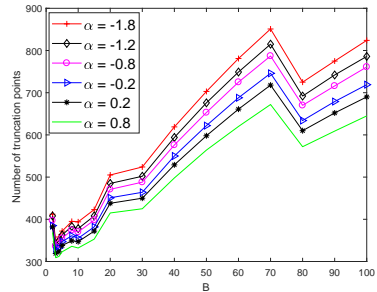,width=5.5cm}
{\par (b)  }
\end{minipage}
\end{center}
\caption{(a) The maximum relative error $\|e\|_{\infty}$ against the basis $B$;
(b) The total number $\sum {q}_{\alpha}(N_{\ell})$ of the  quadrature points
against $B$; $\tau = 0.01,T=10^4,\epsilon=10^{-10}$.
 \label{Fig-2}}
\end{figure}

Figure \ref{Fig-2} (a) shows the relative maximum error $\|e\|_{\infty}$ against
the basis $B$ for different fractional orders,  $\tau = 0.01,T=10^4,\epsilon=10^{-10}$.
We can see that better results are obtained than the predicted precision $\epsilon=10^{-10}$,
even for $\alpha=-1.8$ (the fractional derivative of order $1.8$). Figure \ref{Fig-2} (b)
shows the  total number of the truncated LG quadrature points against the basis $B$.
For a fixed $B$, the number of  truncated quadrature points increases as $\alpha$
decreases, which is in line with the theoretical prediction \eqref{kappa-N}.
For a fixed fractional order $\alpha$, the  number of the truncated quadrature points increases
as $B$ increases, which agrees with \eqref{eq:sec-5-3-5}, due to
$\mathcal{T}_{\ell}\approx 2B-1$ when $\ell$ is sufficiently large.

It is reasonable to choose a relatively smaller
precision $\epsilon$ and smaller basis $B$ in numerical simulations,
since a smaller $\epsilon$ ensures a relatively large  $N_{\ell}$ that guarantees
the exponential convergence of the  LG quadrature, and a small $B$  ensures a
smaller radius of convergence of the quadrature error \eqref{eq:sec-5-3-6}.

%

\begin{table}[!h]
\caption{The maximum relative error
$\|e\|_{\infty}$ under different precision $\epsilon$,  $B= 5$,
$\alpha=-0.5,{\Delta}T = 1$.}\label{tb2-5}
\centering\footnotesize
 \begin{tabular}{|c|ccc|ccc|}
\hline
\multicolumn{1}{|c|}{ } &
\multicolumn{3}{|c|}{$\tau=0.01,T=10^4$}&\multicolumn{3}{|c|}{$\tau=0.1,T=10^5$} \\\hline
$\epsilon$&$\sum N_{\ell}$ &$\sum {q}_{\alpha}(N_{\ell})$ & $\|e\|_{\infty}$
 &$\sum N_{\ell}$ &$\sum {q}_{\alpha}(N_{\ell})$ & $\|e\|_{\infty}$\\\hline
$10^{-12}$ &1018    &384  &1.1937e-12&1203    &442  &6.8208e-12    \\
$10^{-10}$ &848     &349  &4.7044e-12&1001    &402  &6.8191e-12    \\
$10^{-8}$  &680     &313  &4.7044e-12&799     &361  &6.5369e-12    \\
$10^{-6}$  &512     &271  &3.1858e-09&604     &312  &1.2450e-10     \\
$10^{-5}$  &427     &243  &2.1577e-06&503     &281  &1.8889e-08     \\
$10^{-4}$  &340     &215  &2.1577e-06&404     &252  &2.3317e-07     \\
\hline
\end{tabular}
\end{table}

\begin{figure}[!h]
\begin{center}
\begin{minipage}{0.47\textwidth}\centering
\epsfig{figure=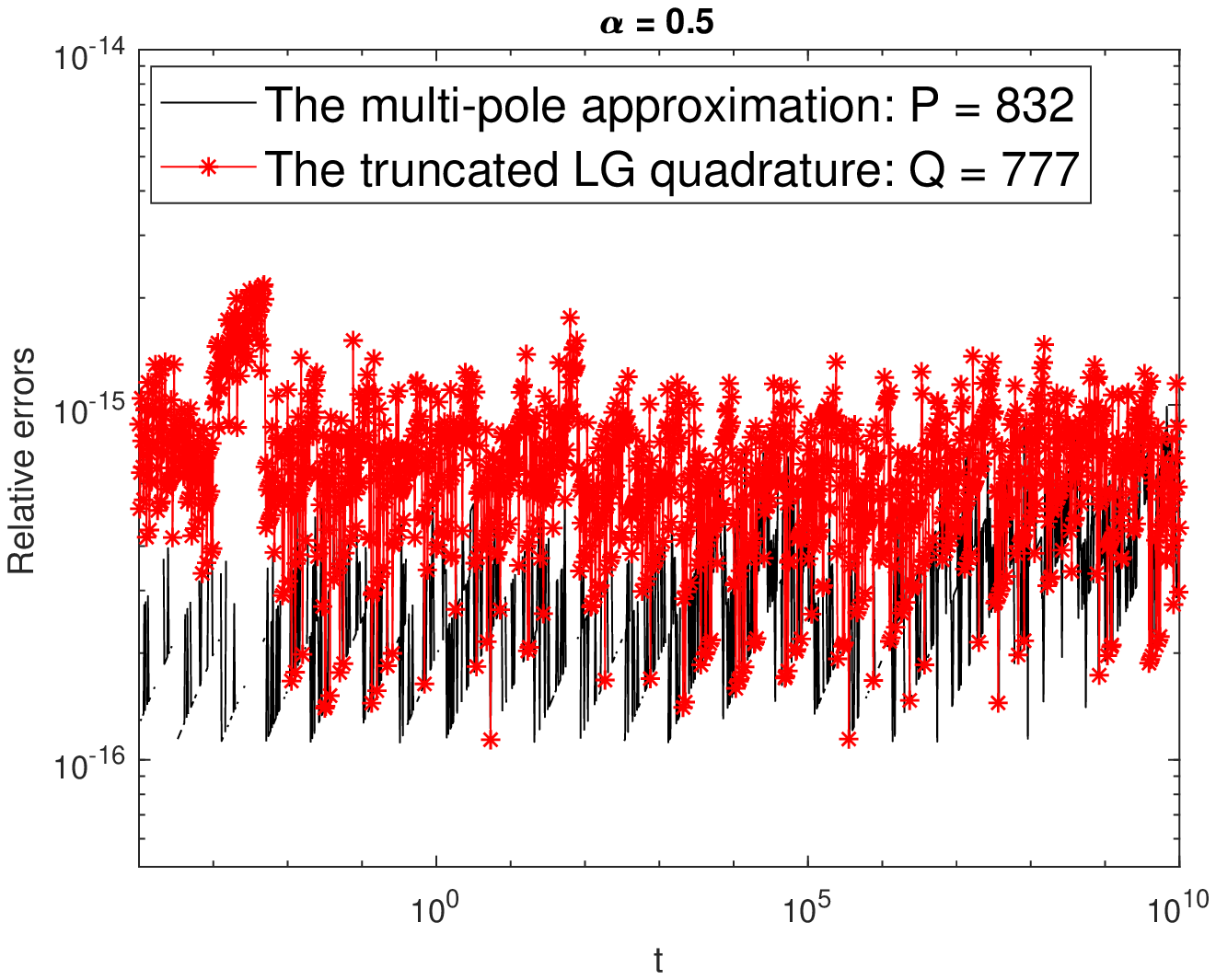,width=5.5cm}
\par {(a)  $\alpha = 0.5,\tau=0.0001$.}
\end{minipage}
\begin{minipage}{0.47\textwidth}\centering
\epsfig{figure=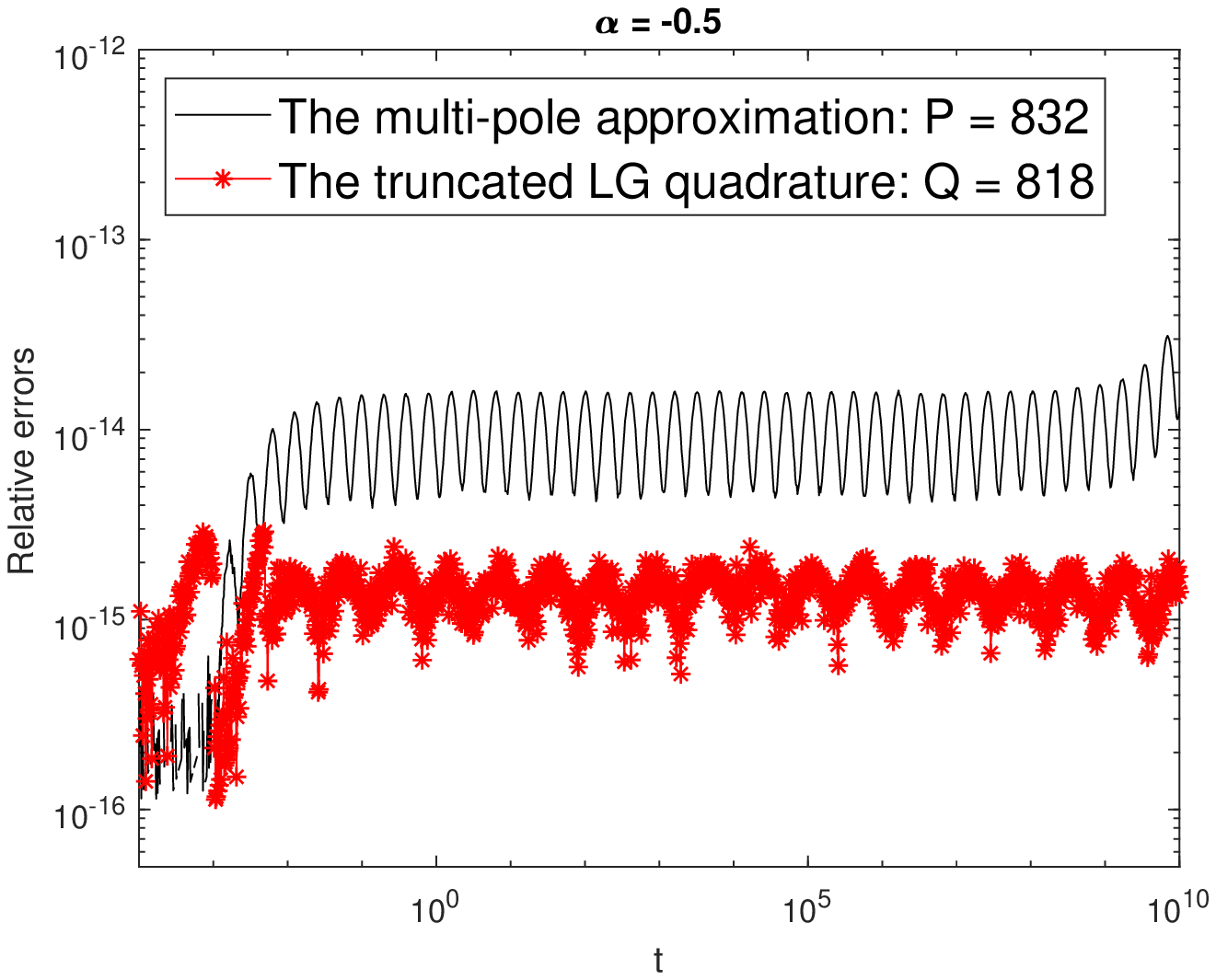,width=5.5cm}
\par {(b)   $\alpha = -0.5,\tau=0.0001$.}
\end{minipage}
\end{center}
\caption{Comparison of the truncated LG quadrature and the multipole approximation
\cite{BafHes17b}. \label{eg32fig1-2}}
\end{figure}

Finally in this section, we compare the truncated LG quadrature with the
multipole approximation proposed in \cite{BafHes17b}.
As the key idea of the existing fast methods aforementioned
is to seek a sum-of-exponentials of the form  $\sum \omega_je^{-\lambda_jt}$
to approximate the kernel function $k_{\alpha}(t)$,
we compare the accuracy of the sum-of-exponentials from the LG quadrature
and the multipole approximation. The relative pointwise errors
$e_n=\big|k_{\alpha}(t_n)-\sum \omega_je^{-\lambda_jt_n}\big|/\big|k_{\alpha}(t_n)\big|$
are shown in Figure \ref{eg32fig1-2}, where we set the precision $\epsilon=10^{-10}$,
$B=5$, and $\Delta T = \tau$ when the LG quadrature is applied,
and the precision  in \cite{BafHes17b}
is set to be $10^{-14}$. For  $\alpha=0.5$, the two methods achieve  similar
accuracy with, respectively,  $P=832$ and $Q=777$ quadrature points
for the multipole approximation and the truncated LG quadrature,
see Figure \ref{eg32fig1-2} (a). Figure \ref{eg32fig1-2} (b) shows
the pointwise errors   for $\alpha = -0.5$, the truncated LG
quadrature shows a slightly better approximation. For other fractional orders
$\alpha \in (-1,1)$, the truncated LG quadrature is competitive with
the multipole approximation both in accuracy and the computational cost. These
results are not shown here.

\section{Numerical examples and applications}\label{sec:numerical}
In this section, two examples are presented to verify the effectiveness
of the present fast method when it is applied to solve FDEs.
All the algorithms are implemented using  MATLAB 2016a, which were
run in a 3.40 GHz PC having 16GB RAM and Windows 7 operating system.
\begin{example}\label{s5-eg-1}
Consider the following scaler FDE
\begin{equation}\label{sec5:eq-1}
{}_{C}D^{\alpha}_{0,t}u(t) = -A u(t) + F(u,t),{\quad}u(0)=u_0,{\quad}t\in(0, T],
\end{equation}
where $0<\alpha \leq 1$ and $A \geq 0$, and ${}_{C}D^{\alpha}_{0,t}$ is the Caputo
fractional operator, which satisfies
${}_{C}D^{\alpha}_{0,t}u(t)=k_{\alpha}*u(t)-{u(0)t^{-\alpha}}/{\Gamma(1-\alpha)}$,
see \cite{Pod-B99}.
\end{example}

We present our fast numerical method for \eqref{sec5:eq-1} as follows:
For a given length $\Delta T=n_0\tau$,  find $U_n$ for $n>n_0$ such that
\begin{equation}\label{sec5:eq-1-2}
{}_FD^{(-\alpha,n,m)}_{\Delta T,\tau}U -  {u_0 t_n^{-\alpha}}/{\Gamma(1-\alpha)}
= -A U_n + F(U_n,t_n),
\end{equation}
where ${}_FD^{(\alpha,n,m)}_{\Delta T,\tau}$ is defined by \eqref{fast-conv-correction}.
The application of the direct convolution method
means here that ${}_FD^{(-\alpha,n,m)}_{\Delta T,\tau}U$
in \eqref{sec5:eq-1-2} is replaced by $D^{(-\alpha,n,m)}_{\Delta T,\tau}U $.
The Newton   method is used to solve the nonlinear system \eqref{sec5:eq-1-2}.
The direct method is applied to obtain $U_k(\leq k\leq m)$, and the
stating values $U_k(1\leq k \leq m)$ are obtained using the direct method
with one correction term and a smaller stepsize $\tau^{-2}$.

The following two cases are considered in this example.
\begin{itemize}
  \item[] Case I: For the linear case of $F=0$, the exact solution of \eqref{sec5:eq-1} is
  $$u(t)=E_{\alpha}(-At^{\alpha}),$$
  where $E_{\alpha}(t)$ is the  Mittag--Leffler function defined by
$E_{\alpha}(t)=\sum_{k=0}^{\infty}\frac{t^{k}}{\Gamma(k\alpha+1)}.$
  \item[] Case II: Let $F=u(1-u^2)$
  and the initial condition is taken as $u_0=1$.
\end{itemize}

The maximum error is defined by
$$\|e\|_{\infty}=\max_{0\leq n \leq {T}/\tau}\big|e_n\big|,  {\quad}e_n=u(t_n)-U_n,T=40.$$

In this example,
we always choose the memory length  ${\Delta T}=0.5$,  the basis $B=5$,
the  truncation number ${q}_{\alpha}{(N_\ell)}$ in \eqref{eq:sec-5-20} is
calculated by \eqref{kappa-N} and $N_\ell$ is determined by \eqref{eq:sec-5-3-5}
under the precision  $\epsilon=10^{-10}$,
and $A=1$. We will reset these parameters if needed.
The fast method based on quadratic interpolation \eqref{quadratic} is applied
in this example  if there is no further illustration.

The purpose of Case I is to check the effectiveness of the present fast method
for non-smooth solutions.
We demonstrate that adding correction terms improves the accuracy of the numerical
solutions significantly.
We first let $\alpha=0.8$, the maximum error and the error at $t=40$ are shown in
Tables  \ref{s5:tb1} and \ref{s5:tb2}, respectively. We can see that the expected
convergence rate $3-\alpha$ is almost achieved by  adding two or three
correction terms,  and the numerical solution at the final time
is much more accurate than that
near the origin.
This phenomenon can be observed from the existing time-stepping methods for
time-fractional differential equations.

As the fractional order decreases, the singularity of the analytical solution becomes stronger.
Tables  \ref{s5:tb3} and \ref{s5:tb4}, respectively, show the maximum errors and the errors
at $t=40$ for $\alpha=0.1$.  We observe that the maximum error, which occurs near the origin,
almost does not improve, though the step size decreases.
We find that as the number of   correction terms increases, the accuracy of the
numerical solutions increases significantly. We need almost 30
correction terms to derive the expected global convergence rate $3-\alpha$ for $\alpha=0.1$,
which is impossible to achieve by   double precision.
The fact is a small number of correction terms is enough to yield satisfactory numerical
solutions, which makes the present method more practical.
We refer readers to \cite{ZengZK17} for more
explanations and numerical results associated with this phenomenon.

\begin{table}[!h]
\caption{{The maximum error $\|e\|_{\infty}$ of the method \eqref{sec5:eq-1-2}, Case I, $\sigma_k=k\alpha$,  $\alpha=0.8$,  ${T=40}$, $B=5$.}}\label{s5:tb1}
\centering\footnotesize
\begin{tabular}{|c|c|c|c|c|c|c|c|c|c|c|c|c|}
\hline
 $\tau$ & $m=0$ & Order& $m=1$ & Order& $m=2$ & Order & $m=3$   & Order\\
 \hline
$2^{-5}$ &2.7966e-3&      &1.5674e-3&      &3.0233e-5&      &3.4818e-5&       \\
$2^{-6}$ &1.7545e-3&0.6726&5.4427e-4&1.5260&5.9024e-6&2.3568&7.7701e-6&2.1638 \\
$2^{-7}$ &1.0604e-3&0.7264&1.8498e-4&1.5569&1.6683e-6&1.8229&1.7199e-6&2.1756 \\
$2^{-8}$ &6.2689e-4&0.7584&6.2086e-5&1.5751&4.2485e-7&1.9733&3.7854e-7&2.1838 \\
$2^{-9}$ &3.6575e-4&0.7773&2.0686e-5&1.5856&1.0205e-7&2.0576&8.3007e-8&2.1892 \\
\hline
\end{tabular}
\end{table}

\begin{table}[!h]
\caption{{The absolute error $\left|e_n\right|$ of the
method \eqref{sec5:eq-1-2}
at ${t=40}$, Case I, $\sigma_k=k\alpha$,  $\alpha=0.8$, $B=5$.}}\label{s5:tb2}
\centering\footnotesize
\begin{tabular}{|c|c|c|c|c|c|c|c|c|c|c|c|c|}
\hline
 $\tau$ & $m=0$ & Order& $m=1$ & Order& $m=2$ & Order & $m=3$   & Order\\
 \hline
$2^{-5}$&5.8890e-7&      &1.9977e-7&      &1.7531e-7 &      &4.8407e-6&      \\
$2^{-6}$&3.0861e-7&0.9322&5.8512e-8&1.7715&3.8771e-8 &2.1769&1.1237e-6&2.1069\\
$2^{-7}$&1.5741e-7&0.9713&1.6977e-8&1.7852&8.5112e-9 &2.1875&2.5398e-7&2.1455\\
$2^{-8}$&7.9372e-8&0.9878&4.8985e-9&1.7931&1.8745e-9 &2.1829&5.6424e-8&2.1703\\
$2^{-9}$&3.9821e-8&0.9951&1.4075e-9&1.7992&4.3619e-10&2.1035&1.2218e-8&2.2073\\
\hline
\end{tabular}
\end{table}

\begin{table}[!h]
\caption{{The maximum error $\|e\|_{\infty}$ of the method \eqref{sec5:eq-1-2}, Case I, $\sigma_k=k\alpha$,  $\alpha=0.1$,  ${{T}=40}$, $B=5$.}}\label{s5:tb3}
\centering\footnotesize
\begin{tabular}{|c|c|c|c|c|c|c|c|c|c|c|c|c|}
\hline
 $\tau$ & $m=0$ & Order& $m=1$ & Order& $m=3$ & Order & $m=5$   & Order\\
 \hline
$2^{-5}$&2.5872e-3&      &1.5852e-3&      &2.1191e-5&      &2.1191e-5&      \\
$2^{-6}$&2.5323e-3&0.0310&1.5017e-3&0.0781&9.6408e-6&1.1362&9.6408e-6&1.1362\\
$2^{-7}$&2.4750e-3&0.0330&1.4177e-3&0.0831&4.3479e-6&1.1488&4.3479e-6&1.1488\\
$2^{-8}$&2.4148e-3&0.0355&1.3338e-3&0.0880&1.9444e-6&1.1610&1.9444e-6&1.1610\\
$2^{-9}$&2.3516e-3&0.0383&1.2507e-3&0.0928&1.2279e-6&0.6631&8.6285e-7&1.1721\\
\hline
\end{tabular}
\end{table}

\begin{table}[!h]
\caption{{The absolute  error $\left|e_n\right|$ of the method \eqref{sec5:eq-1-2}
at  ${t=40}$, Case I, $\sigma_k=k\alpha$,  $\alpha=0.1$, $B=5$.}}\label{s5:tb4}
\centering\footnotesize
\begin{tabular}{|c|c|c|c|c|c|c|c|c|c|c|c|c|}
\hline
 $\tau$ & $m=0$ & Order& $m=1$ & Order& $m=3$ & Order & $m=5$   & Order\\
 \hline
$2^{-5}$&6.4561e-6&      &6.6140e-7&      &1.0065e-8 &      &4.1985e-9 &      \\
$2^{-6}$&3.2045e-6&1.0106&3.1760e-7&1.0583&3.8428e-9 &1.3891&1.0230e-9 &2.0371\\
$2^{-7}$&1.5919e-6&1.0094&1.5234e-7&1.0599&1.5471e-9 &1.3126&2.5282e-10&2.0166\\
$2^{-8}$&7.9116e-7&1.0087&7.2987e-8&1.0616&6.4397e-10&1.2645&6.4081e-11&1.9801\\
$2^{-9}$&3.9331e-7&1.0083&3.4930e-8&1.0632&2.7313e-10&1.2374&1.6914e-11&1.9217\\
\hline
\end{tabular}
\end{table}

Table \ref{s5:tb4-2} shows the difference
$\eta=\max_{0\leq n \leq T/\tau}{\big|U_D^n-{U}_F^n\big|}$,
where  $U^n_F$   are numerical solutions derived from
the fast method under the precision $\epsilon$,
 and $U^n_D$ are numerical solutions  from the direct method.
 We can see that the difference $\eta$ is
independent of the stepsize $\tau$, which confirms \eqref{eq:error-2}.
Table \ref{s5:tb4-2} also shows that the accuracy of the present fast calculation
outperforms the predicted accuracy $\epsilon$, see \eqref{eq:error-2}.

\begin{table}[!h]
\caption{{The  difference  of
the numerical solutions between
the fast method and the direct method, Case I,
$m=0$,  $\alpha=0.1$,  ${{T}=40}$, $B=5$.}}\label{s5:tb4-2}
\centering\footnotesize
\begin{tabular}{|c|c|c|c|c|c|c|c|c|c|c|c|c|}
\hline
$\epsilon$& $\tau=2^{-5}$ & $\tau=2^{-6}$ & $\tau=2^{-7}$& $\tau=2^{-8}$& $\tau=2^{-9}$ \\
 \hline
$10^{-12}$&2.8255e-13&2.7850e-13&2.7167e-13&1.1102e-15&1.4988e-15\\
$10^{-10}$&2.8239e-13&2.7839e-13&2.7162e-13&5.6066e-15&2.1094e-15\\
$10^{-8}$ &2.8116e-13&3.8219e-13&6.6397e-13&2.2093e-14&8.2682e-13\\
$10^{-6}$ &7.5677e-11&7.2366e-11&8.0116e-10&2.6198e-11&8.0352e-12\\
$10^{-5}$ &3.2916e-11&1.9709e-10&1.2477e-10&2.6461e-10&2.2769e-10\\
$10^{-4}$ &3.0868e-09&9.4674e-09&2.4896e-09&1.5458e-08&4.3178e-08\\
\hline
\end{tabular}
\end{table}

In Table \ref{s5:tb5}, we compare the present fast method \eqref{sec5:eq-1-2}
based on the linear interpolation with the graded mesh method in  \cite{StynesORiGra17},
in which the nonuniform grid points are given by $t_j=(j\tau)^r$,
$0\leq j \leq 1/\tau$. We can see that the fast method with  correction terms
is competitive with the graded mesh method for $\alpha=0.5$. In \cite{StynesORiGra17},
the authors obtained the optimal grid mesh  $t_j=(j\tau)^{(2-\alpha)/\alpha}$
to achieve the global $(2-\alpha)$-order accuracy, which works well
when $\alpha$ is relatively large and the expected convergence rate is achieved.
For the correction method, too many correction terms may harm the accuracy,
but a few number of correction terms can achieve highly accurate numerical
solutions, which is not investigated here, readers can refer to
\cite{DieFFW06,Lub86,ZengZK17} for more discussion.
Note that for values of $\alpha$ close to zero, the correction method is even more
effective than the optimal graded mesh approach in \cite{StynesORiGra17}.


\begin{table}[!h]
\caption{{Comparison of the present method
with the graded mesh method in \cite{StynesORiGra17} based on linear interpolation, Case I,
$\sigma_k=k\alpha$,  $\alpha=0.5$, $B=5$. The errors $\|e\|_{\infty}$
are displayed.}}\label{s5:tb5}
\centering\footnotesize
The graded mesh method \cite{StynesORiGra17} with grid points $t_j=(j\tau)^r$.
\begin{tabular}{|c|c|c|c|c|c|c|c|c|c|c|c|c|}
\hline
 $\tau$ & $r=1$ & Order& $r=3/2$ & Order& $r=3$ & Order & $r=6$   & Order\\
 \hline
$2^{-5}$&3.6491e-2&      &1.6839e-2&      &2.6568e-3&      &2.5431e-3&       \\
$2^{-6}$&2.7048e-2&0.4320&1.0288e-2&0.7108&1.0077e-3&1.3986&9.2945e-4&1.4521 \\
$2^{-7}$&1.9772e-2&0.4521&6.2162e-3&0.7268&3.7277e-4&1.4347&3.3644e-4&1.4660 \\
$2^{-8}$&1.4312e-2&0.4663&3.7315e-3&0.7363&1.3582e-4&1.4566&1.2096e-4&1.4759 \\
$2^{-9}$&1.0288e-2&0.4762&2.2313e-3&0.7419&4.9041e-5&1.4696&4.3280e-5&1.4827 \\
\hline
\end{tabular}
\centering\footnotesize
{
The fast method based on linear interpolation with correction terms
\begin{tabular}{|c|c|c|c|c|c|c|c|c|c|c|c|c|}
\hline
 $\tau$ & $m=1$ & Order& $m=2$ & Order& $m=3$ & Order & $m=4$   & Order\\
 \hline
$2^{-5}$&5.6366e-3&      &2.1893e-4&      &1.4976e-4&      &7.2351e-5&       \\
$2^{-6}$&3.1659e-3&0.8322&1.0737e-4&1.0279&6.9115e-5&1.1156&4.4890e-5&0.6886 \\
$2^{-7}$&1.7231e-3&0.8777&4.9196e-5&1.1260&2.9276e-5&1.2393&2.2604e-5&0.9898 \\
$2^{-8}$&9.1600e-4&0.9116&2.1427e-5&1.1991&1.1719e-5&1.3208&1.0105e-5&1.1614 \\
$2^{-9}$&4.7862e-4&0.9364&8.9763e-6&1.2552&4.5167e-6&1.3756&4.1916e-6&1.2696 \\
\hline
\end{tabular}}
\end{table}

\begin{figure}[!h]
\begin{center}
\begin{minipage}{0.47\textwidth}\centering
\epsfig{figure=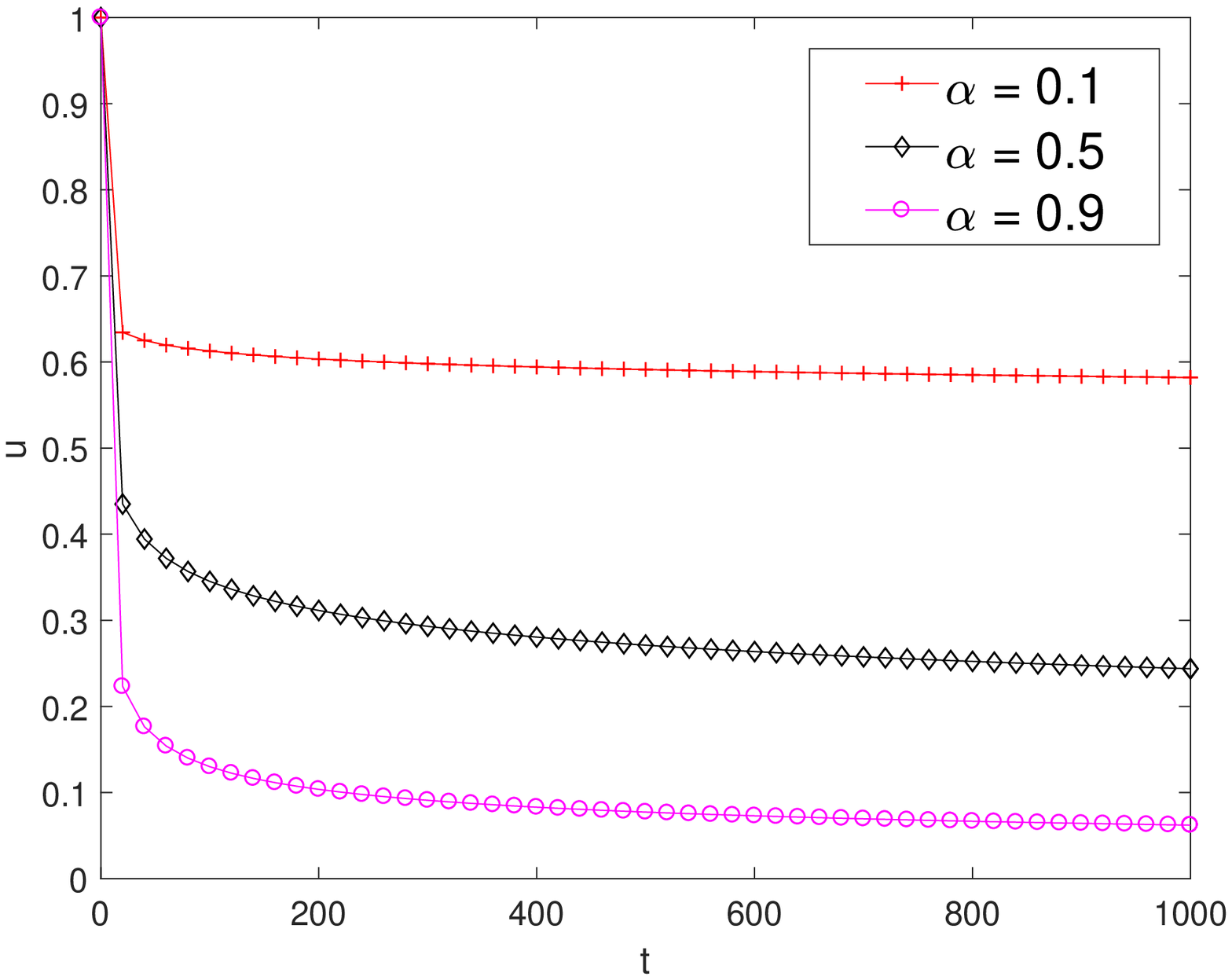,width=5.5cm}
\par {(a) Numerical solutions. }
\end{minipage}
\begin{minipage}{0.47\textwidth}\centering
\epsfig{figure=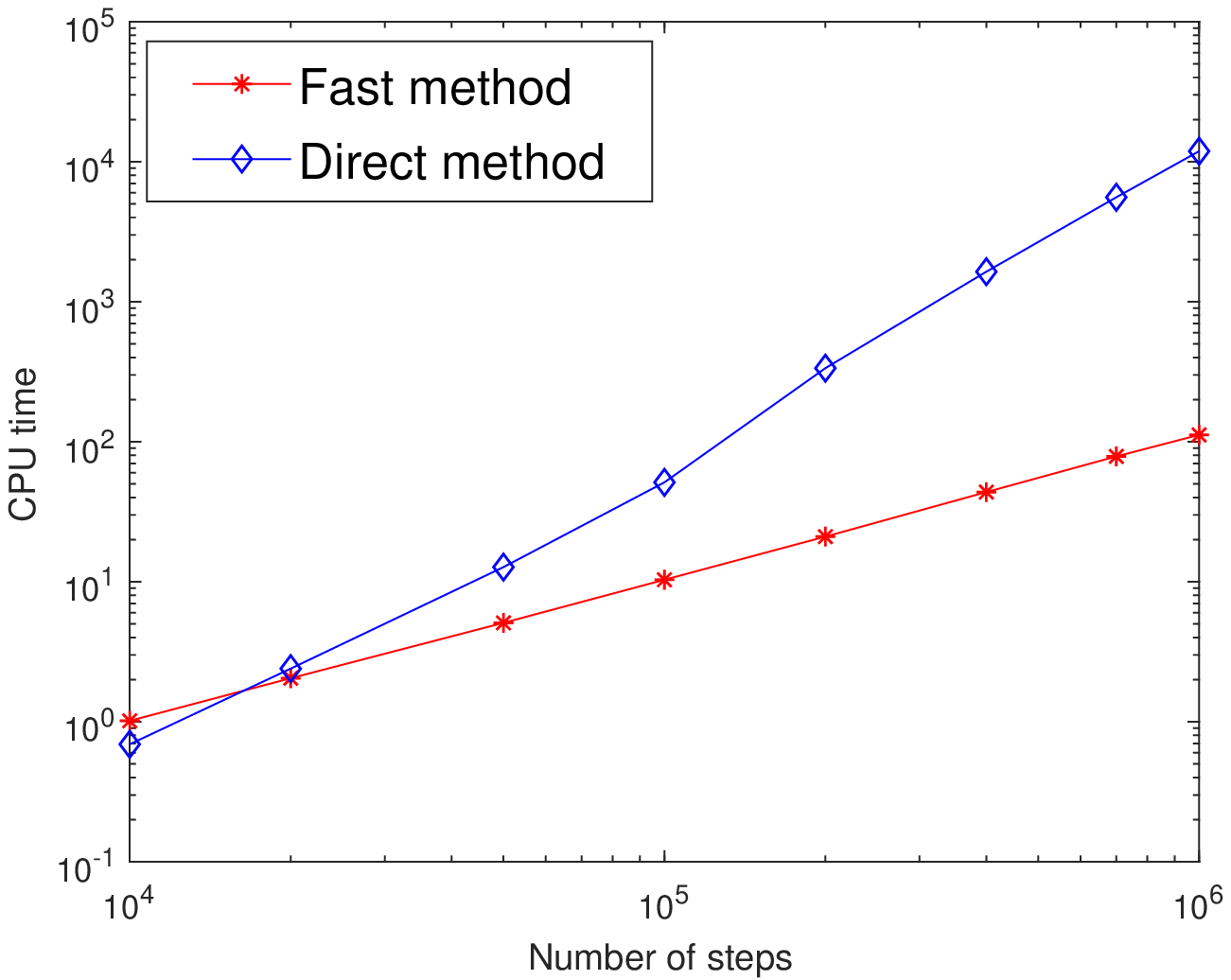,width=5.5cm}
\par {(b)   Computational time.}
\end{minipage}
\end{center}
\caption{Numerical solutions and the computational time of
the two different convolutions for Case II,
$\tau=0.01,B=5,m=2,\sigma_k=k\alpha$.\label{eg51fig1}}
\end{figure}

Next, we implement the present fast method for a longer time computation and
compare it  with the direct convolution method.
We show in Figure \ref{eg51fig1}(a)  numerical solutions for $t\in[0,{T}],{T}=10000$,
where we choose $\alpha=0.1,0.5,0.9$, and the time stepsize $\tau=0.01$.
In Figure \ref{eg51fig1}(b), we plot the computational time of the
fast method and the direct method with two correction terms applied.
It shows that the computational cost of the fast convolution almost achieves linear complexity,
which is much less  than that of the direct convolution with
 $O(n_T^2)$ operations when $n_T$ is sufficiently large,
 where $n_0={\Delta T}/\tau$ and $n_T=T/\tau$.

Next, we consider a system of FODEs each having a different fractional index.
\begin{example}\label{s5-eg-2}
Consider the following   system of fractional differential equations
each with different fractional index
\begin{equation}\label{sec5:eq-2}\left\{\begin{aligned}
&{}_{C}D^{\alpha_1}_{0,t}u(t) = w + (v-c_1)u,\\
&{}_{C}D^{\alpha_2}_{0,t}v(t) =  1- c_2v - u^2,\\
&{}_{C}D^{\alpha_3}_{0,t}w(t) = -u - c_3w,
\end{aligned}\right.\end{equation}
where $0< \alpha_k \leq 1\,(k=1,2,3)$, $c_1,c_2$ and $c_3$ are positive constants, $c_2>1/2$.
\end{example}

Let $U_n,V_n$ and $W_n$ be the approximate solutions
of $u(t_n),v(t_n)$ and $w(t_n)$, respectively. Then the fully discrete scheme  for
the system \eqref{sec5:eq-2} is given by:
\begin{equation}\label{sec5:eq-3}\left\{\begin{aligned}
&{}_FD^{(-\alpha_1,n,m)}_{\Delta T,\tau}U-{u_0 t_n^{-\alpha_1}}/{\Gamma(1-\alpha_1)}
= W_n + (V_n-c_1)U_n,\\
&{}_FD^{(-\alpha_2,n,m)}_{\Delta T,\tau}V - {v_0 t_n^{-\alpha_2}}/{\Gamma(1-\alpha_2)}
 = 1- c_2V_n - U_n^2,\\
&{}_FD^{(-\alpha_3,n,m)}_{\Delta T,\tau}W - {w_0 t_n^{-\alpha_3}}/{\Gamma(1-\alpha_3)}
= -U_n - c_3W_n,
\end{aligned}\right.\end{equation}
where $n\geq 2$ and ${}_FD^{(\alpha,n,m)}_{\Delta T,\tau}U$
is defined by \eqref{fast-conv-correction}.
Because we use  quadratic interpolation,
$U_k,V_k,W_k$ for $k=1,2$ need to be known, which can be derived using the known methods
with smaller stepsize. Here we use the L1 method with  one correction term and smaller
stepsize $\tau^2$ to derive these values.
We  take   ${\Delta T}=\tau$ and $B=5$ in numerical simulations of this example.

\begin{figure}[!h]
\begin{center}
\begin{minipage}{0.47\textwidth}\centering
\epsfig{figure=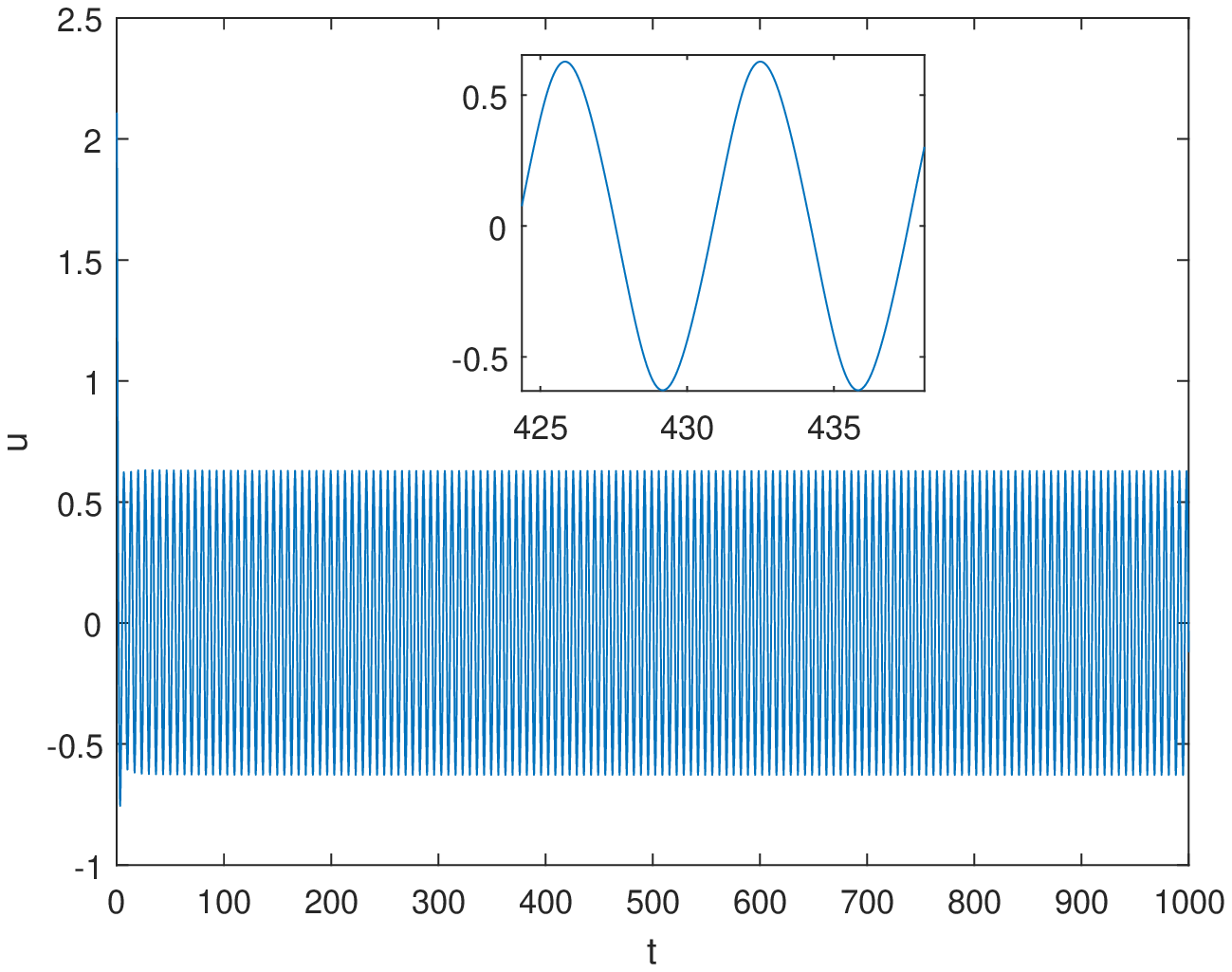,width=5.5cm}
\par {(a)   }
\end{minipage}
\begin{minipage}{0.47\textwidth}\centering
\epsfig{figure=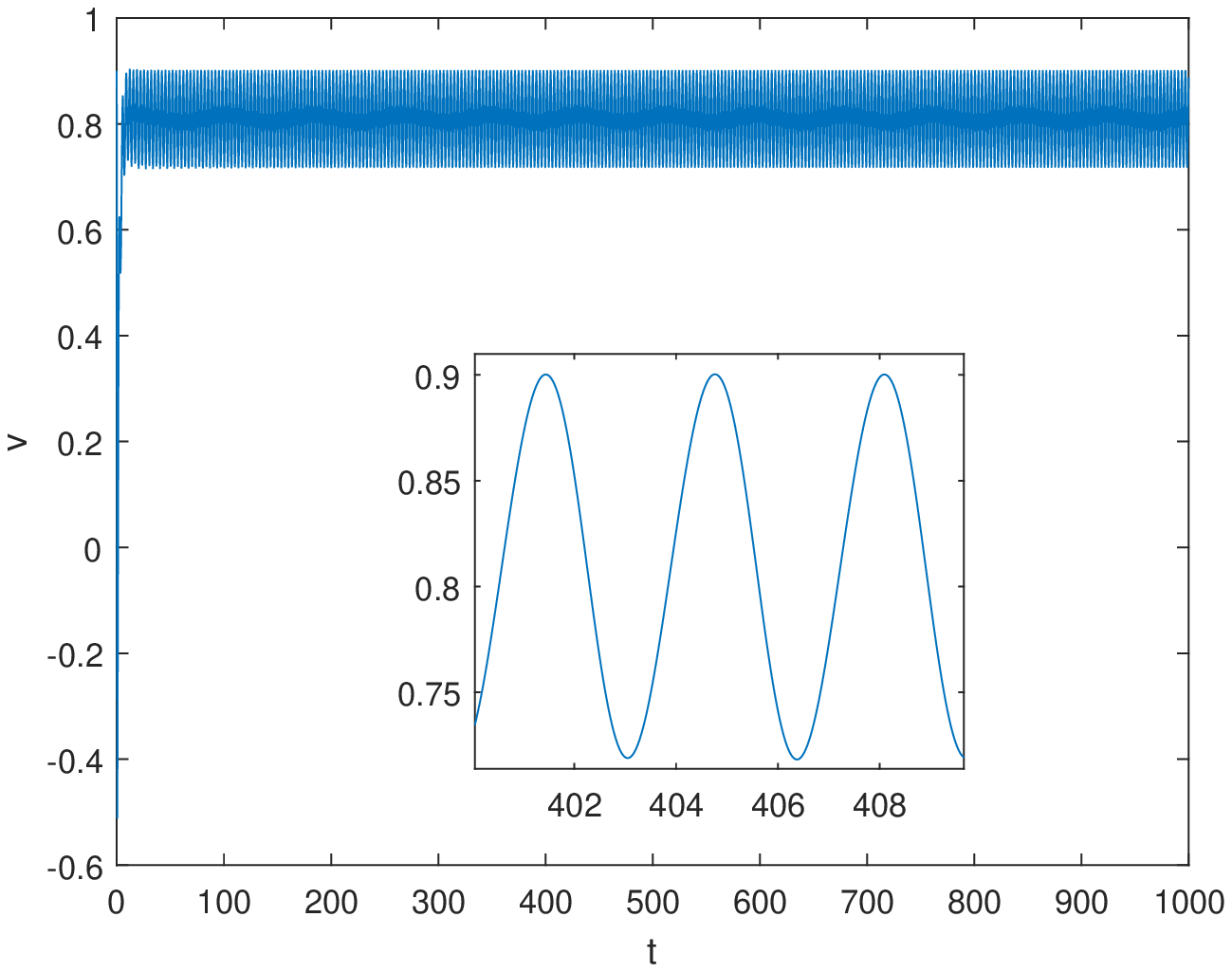,width=5.5cm}
\par {(b)  }
\end{minipage}
\begin{minipage}{0.47\textwidth}\centering
\epsfig{figure=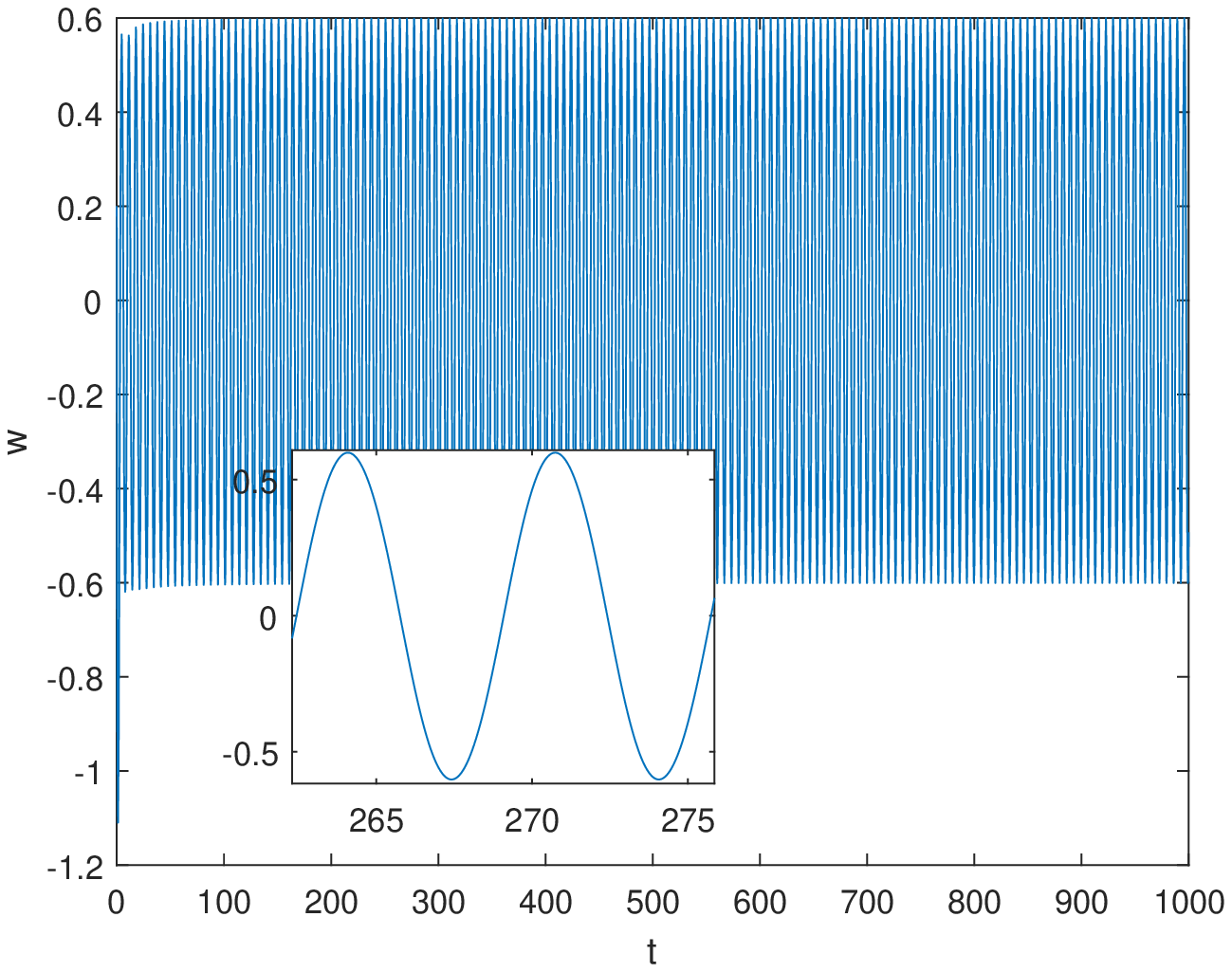,width=5.5cm}
\par {(c)   }
\end{minipage}
\begin{minipage}{0.47\textwidth}\centering
\epsfig{figure=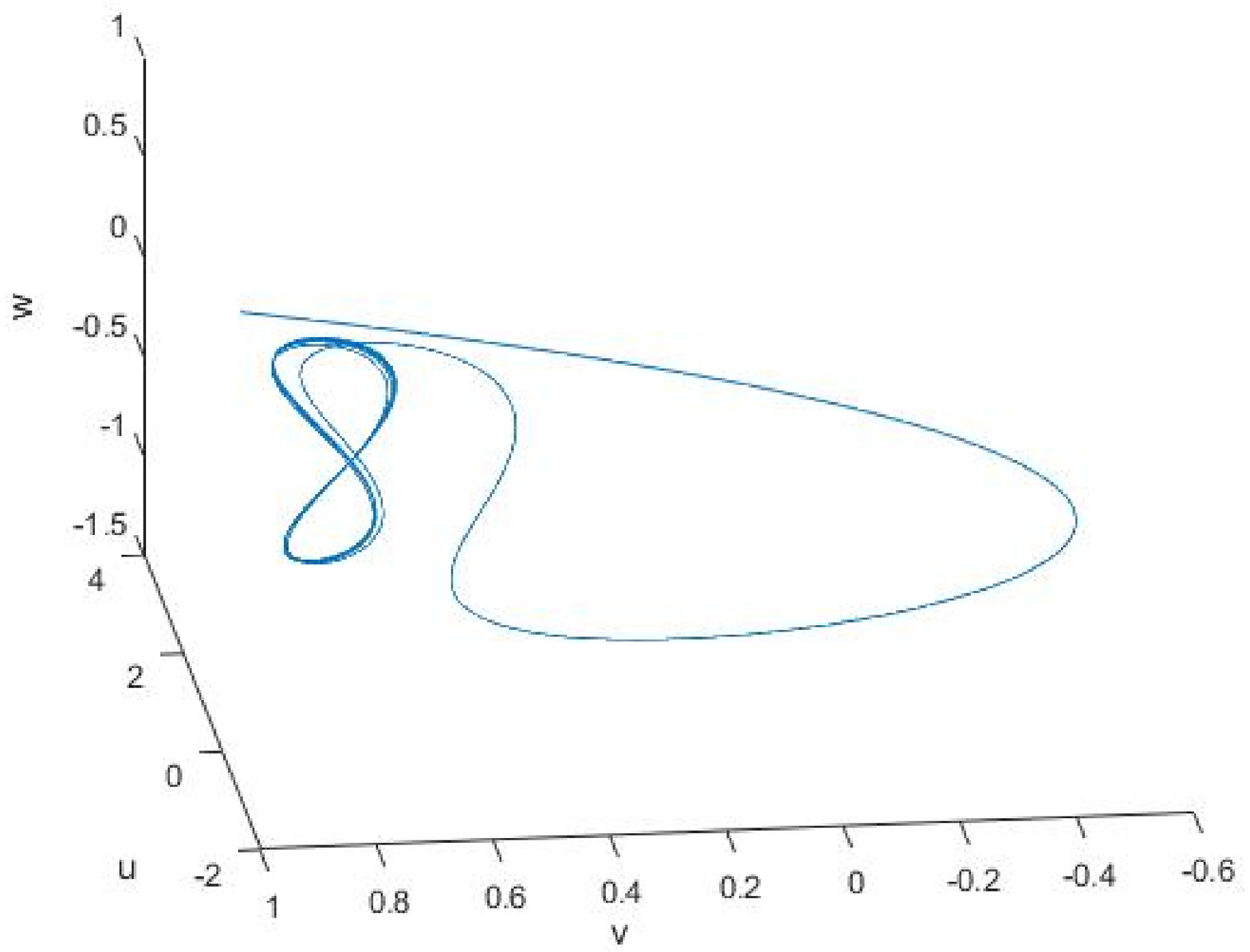,width=5.5cm}
\par {(d)    }
\end{minipage}
\end{center}
\caption{Numerical solutions for Example \ref{s5-eg-2}, $\tau=0.01,m=2,\sigma_k=k\alpha,\alpha=0.9$.\label{eg62fig01}}
\end{figure}

\begin{figure}[!h]
\begin{center}
\begin{minipage}{0.47\textwidth}\centering
\epsfig{figure=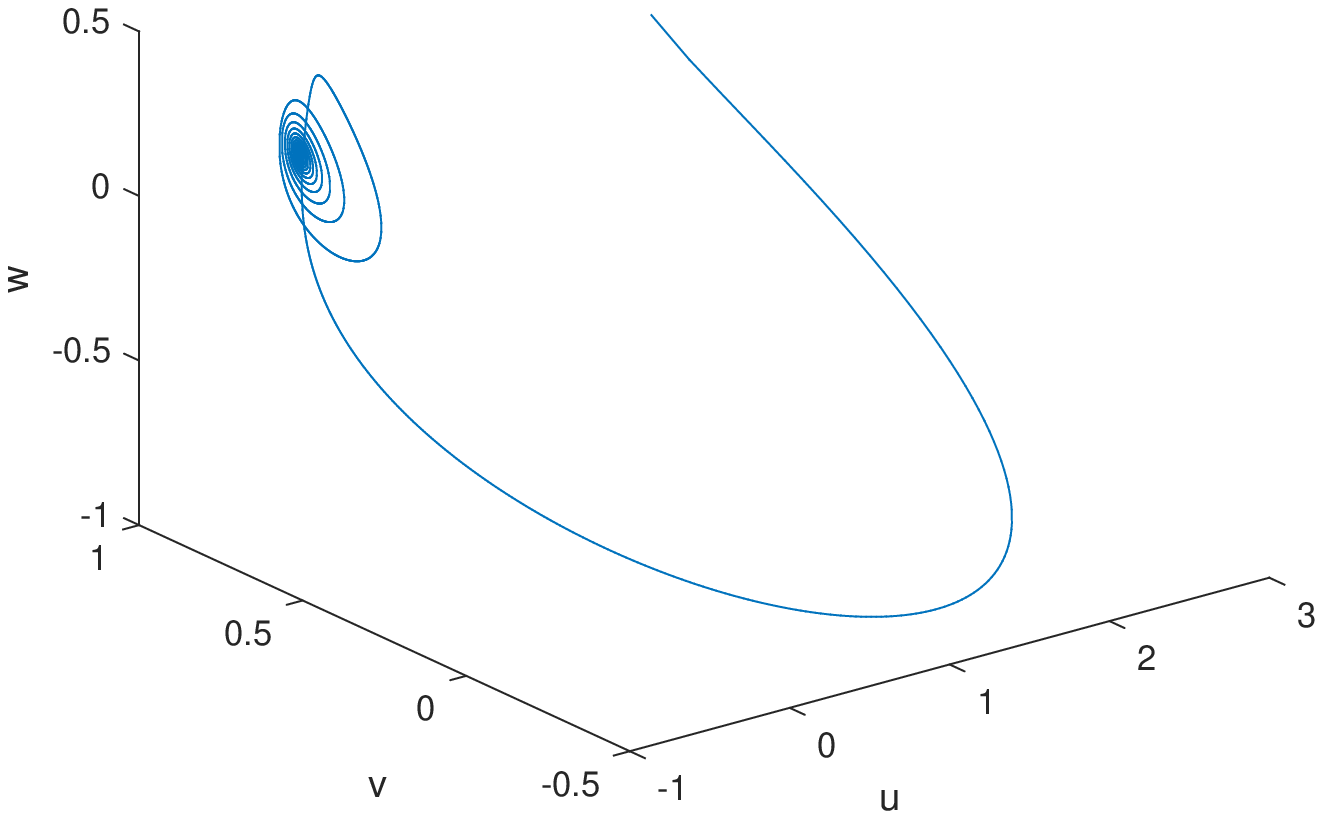,width=5.5cm}
\par {(a) $(\alpha_1,\alpha_2,\alpha_3)=(0.9,0.8,0.7)$. }
\end{minipage}
\begin{minipage}{0.47\textwidth}\centering
\epsfig{figure=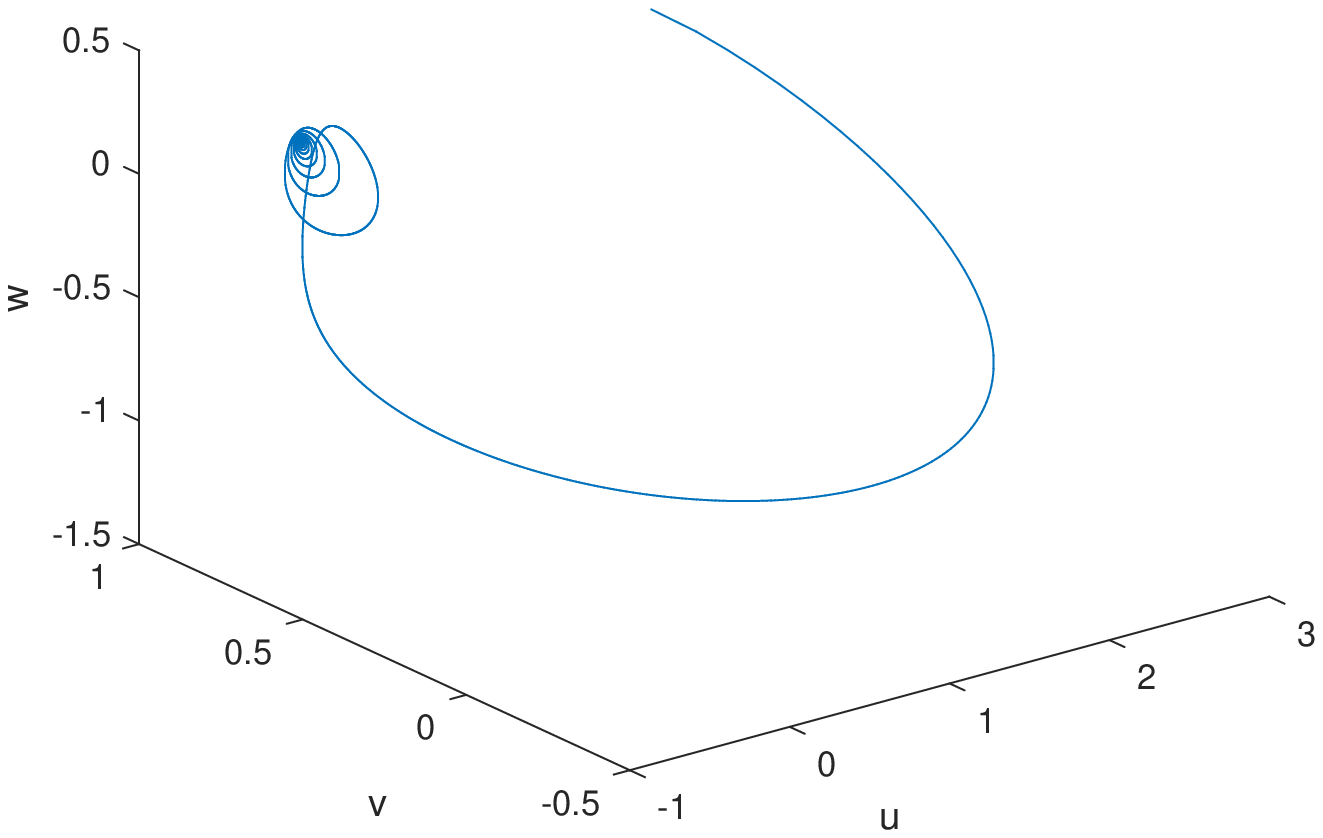,width=5.5cm}
\par {(b)  $(\alpha_1,\alpha_2,\alpha_3)=(0.7,0.8,0.9)$. }
\end{minipage}
\end{center}
\caption{Numerical solutions for Example \ref{s5-eg-2}, $t\in [0,1000]$,
 $\tau=0.01,m=0$.\label{eg62fig03}}
\end{figure}

If $\alpha_1=\alpha_2=\alpha_3=\alpha$, then \eqref{sec5:eq-2}
is just the fractional Lorenz system  \cite{WangXiao2015}.
It  has been proved that the fractional Lorenz system \eqref{sec5:eq-2}
is dissipative \cite{WangXiao2015} and has an absorbing set defined by
a ball $B(0,\sqrt{a/b}+\hat{\epsilon})$, where $a=1/2$ and $b=\min\{c_1,c_2-1/2,c_3\}$.
We take $c_1=1/4,c_2=1,c_3=1/4$ and the
initial conditions as those in \cite{WangXiao2015}, i.e., $U_0=u(0)=2,V_0=v(0)=0.9,W_0=w(0)=0.2$,
the time stepsize is taken as $\tau=0.01$. We compute numerical solutions for
$t\in[0,1000]$, which is much larger than that in  \cite{WangXiao2015}. Numerical solutions
for $\alpha=0.9$ are shown in Figure \ref{eg62fig01}.  We can easily find that
$U_n^2+V_n^2+W_n^2<2$, which means $(U_n,V_n,W_n)\in B(0,\sqrt{2})$.
For other fractional orders $\alpha_k=\alpha$, we have similar results,
see also  \cite{WangXiao2015}. Next, we choose different $\alpha_1,\alpha_2$, and $\alpha_3$,
and exhibit the numerical solutions   in Figure  \ref{eg62fig03}
for $(\alpha_1,\alpha_2,\alpha_3)=(0.9,0.8,0.7)$ and
$(\alpha_1,\alpha_2,\alpha_3)=(0.7,0.8,0.9)$.
We can see that the numerical solutions $(U_n,V_n,W_n)$ are also in a ball.
For other choices of the fractional orders, we have similar results, which
are not provided here.

\section{Conclusion and discussion}\label{concl}
We propose a unified fast memory-saving time-stepping method for
both   fractional integral and derivative operators, which is applied to solve FDEs.
We generalize the fast convolution in \cite{LubSch02} for fractional operators,
while we use the (truncated) LG
quadrature to discretize the kernel in the fractional operators instead of the
trapezoidal rule in \cite{LubSch02}.  We also introduced correction terms in the fast method
such that the non-smooth solutions of the considered FDEs can be resolved accurately.
The present fast method has  $O(n_0+\sum_{\ell}^L{q}_{\alpha}(N_{\ell}))$  active memory and
$O(n_0n_T+ (n_T-n_0)\sum_{\ell}^L{q}_{\alpha}(N_{\ell}))$ operations.
If a suitable
memory length $\Delta T$ and a very large basis $B$ are applied,
then the present fast method is similar to that
in \cite{JiangZZZ17,JingLi10},   that is,  the kernel $k_{\alpha}(t)$ is approximated
via the truncated LG quadrature for all  $t \in [\Delta T, T]$.

We present the error analysis of the present fast method. For a given precision,
a criteria on how to choose the parameters used in the fast method is given explicitly,
which works very well in numerical simulations.

We considered the fast method for the fractional  operator of order $\alpha<1$.
In fact, our method can be extended to
the fractional integral of order greater than one. For example, for $\alpha\in [1,2)$,
Eq. \eqref{kernel} still holds in the sense of the finite part integral, that is,
Eq. \eqref{kernel} is equivalent to $k_{\alpha}(t)=\frac{t}{\Gamma(\alpha)\Gamma(2-\alpha)}
\int_0^{\infty}\lambda^{1-\alpha}e^{-t\lambda}\dx[\lambda].$
In such a case,
a coupled ODEs are needed to be resolved instead of one for $\alpha<1$.

In the future, we will apply the fast method to solve
a large system of time-fractional PDEs, perhaps in three spatial dimensions.
We will also explore how to efficiently calculate
the discrete convolution  $\sum_{k=0}^n\omega_{n-k}u(t_k)$, where the quadrature weights
$\{\omega_{n}\}$ are not from the interpolation as done in the present work,
but from the generating functions,  see,  \cite{Lub86}.

\appendix
\section{Proofs}\label{appendix-A}
Proof of Theorem \ref{thm-1}.
\begin{proof}
We follow the proof of Theorem 2.2 in \cite{Xiang2012}.
We first expand $g(\lambda)=e^{-t\lambda}$ in terms of the Laguerre polynomials, i.e.,
$g(\lambda)=\sum_{n=0}^{\infty}a_nL_n^{(\alpha)}(\lambda),$
where
$$a_n=\frac{\Gamma(n+1)}{\Gamma(n+1+\alpha)}
\int_0^{\infty}\lambda^{\alpha} e^{-\lambda} g(\lambda)\dx[\lambda]
=\frac{t^n}{(t+1)^{n+1+\alpha}}.$$
The following property will be used, see \cite{Xiang2012},
\begin{equation*}
\big|Q_N^{\alpha}[L_n^{(\alpha)}]\big|
\leq \left\{\begin{aligned}
&2\Gamma(1+\alpha),{\quad}-1<\alpha \leq 0,\\
&2^{1+\alpha}{\Gamma(n+1+\alpha)}/{\Gamma(n+1)},{\quad}\alpha>0.
\end{aligned}\right.\end{equation*}
With  the above two equations and
$I^{\alpha}[e^{-t\lambda}]- Q_N^{\alpha}[e^{-t\lambda}]
=\sum_{n=2N}^{\infty}a_nQ_N^{\alpha}[L_n^{(\alpha)}]$ gives
\begin{equation}\label{eq:A-1}
\big|I^{\alpha}[e^{-t\lambda}]- Q_N^{\alpha}[e^{-t\lambda}]\big|
\leq  \left\{\begin{aligned}
&2^{1+\alpha}\Gamma(1+\alpha)\sum_{n=2N}^{\infty}a_n ,{\quad}-1<\alpha \leq 0,\\
&c_{\alpha}2^{1+\alpha}\sum_{n=2N}^{\infty}n^{\alpha}a_n,{\quad}\alpha>0.
\end{aligned}\right.\end{equation}
Let $q=t/(1+t)$. Then we have $\sum_{n=2N}^{\infty}a_n=(1+t)^{-\alpha}q^{2N}$ for $-1<\alpha\leq 0$.
For $\alpha>0$, one hase
$$\sum_{n=2N}^{\infty}n^{\alpha}a_n=(1+t)^{-\alpha}q^{2N}(2N)^{\alpha}
\sum_{n=0}^{\infty}q^n\Big(1+\frac{n}{2N}\Big)^{\alpha}
\leq C_{\alpha,t} 2^{\alpha}(1+t)^{1-\alpha}q^{2N}N^{\alpha},$$
where $C_{\alpha,t} =\sum_{n=0}^{\infty}(n+2)^{\alpha}q^{2Nn}$ is used.
With the above equation and \eqref{eq:A-1} yields \eqref{eq:sec-5-2-2} for $T=1$.
Using the following relation
$$\big|I^{\alpha}[T,e^{-t\lambda}]- Q_N^{\alpha}[T,e^{-t\lambda}]\big|
=T^{-\alpha-1}\big|I^{\alpha}[e^{-(t/T)\lambda}]- Q_N^{\alpha}[e^{-(t/T)\lambda}]\big| $$
leads to \eqref{eq:sec-5-2-2}  for any $T>0$.
The proof is complete.
\end{proof}

Proof of Theorem \ref{thm-2}.
\begin{proof}
The following results can be found in \cite{Gatteschi2012},
\begin{eqnarray}
&&{\lambda}_j < \frac{2j+\alpha+3}{2N+\alpha+3}
\Big(2j+\alpha+3+\sqrt{(2j+\alpha+3)^2 +0.25-\alpha^2}\Big),\label{Laguerre-2}\\
&& {\lambda}_j >  \frac{2(J_{\alpha,j}/2)^2}{2N+\alpha+3},{\quad}
J_{\alpha,j}=\pi(j+3/4+\alpha/2)+O(j^{-1}),\label{Laguerre-3}
\end{eqnarray}
where \eqref{Laguerre-3} holds when $j$ is sufficiently large.
For a sufficiently large $N$,   the Laguerre polynomial satisfies
(see (7.14) in \cite{ShenTW-B11})
\begin{equation}\label{Laguerre-4}
\big|L^{(\alpha)}_N(x)\big| \approx \pi^{-1/2}(Nx)^{-1/4}e^{x/2}, {\quad}\forall x\geq0.
\end{equation}

With \eqref{Laguerre-2}--\eqref{Laguerre-4} and \eqref{Laguerre-weghts}, we have
the following estimate
\begin{equation}\label{Laguerre-5}
{w}^{(\alpha)}_j\leq C_{\alpha}(N+1)^{\alpha} \left(\frac{j+1}{N+1}\right)^3e^{-{\lambda}_j}
\leq C(N+1)^{\alpha} e^{-{\lambda}_j},
\end{equation}
for sufficiently large $j$, where ${\lambda}_j=\theta_j(j+1)^2/(N+1)$   and
$\theta_j$ is bounded and approximately between $\pi^2/4$ and $4$.
The proof is completed.
\end{proof}


\def\cprime{$'$} \def\cprime{$'$} \def\cprime{$'$} \def\cprime{$'$}

\end{document}